\documentclass[twocolumn]{autart}
\usepackage{graphicx}       % Include images
\usepackage{tipa}           % For \textcrb and other phonetic symbols
\usepackage{comment}        % Block comments
\usepackage{xcolor}         % Colored text
\usepackage{lipsum}         % Filler text generation
\usepackage{amsfonts}       % Mathematical fonts
\usepackage[utf8]{inputenc} % Handle strange symbols
\usepackage{amsmath}        % Mathematical symbols
\usepackage{mathtools}      % Add mathematical tricks
\usepackage{amssymb}        % Additional math symbols
\usepackage{alphabeta}      % Greek symbols in text mode
\usepackage{url}            % Handle URLs
\usepackage{nameref}        % Reference section names
\usepackage{cases}          % Cases in equations
\usepackage{xfrac}          % Compact fractions
\usepackage{float}          % Precise figure placement
\usepackage{anyfontsize}    % Flexible font sizes
\usepackage{times}          % More fonts
\usepackage{colortbl}       % Color tables
\usepackage{natbib}         % Necessary bibliography
\usepackage{tikz}           % To make images in latex
\usetikzlibrary{shapes.geometric, arrows.meta, positioning, fit, backgrounds, calc, math, chains}    % Import tikz stuff
\usepackage{pgfmath}
\usepackage{pgfplots}
\usepgfplotslibrary{statistics, groupplots, external, fillbetween}
\usepackage{soul}           % Get strikethrough text
\usepackage[dvipsnames, svgnames]{xcolor}
\pgfplotsset{compat=1.18}
\usepackage{soul}
\usepackage{algorithm}      % Algorithm environments
\usepackage{algpseudocode} 	% Pseudocode formatting (algpseudocodex seems to not work with autart)
\usepackage{etoolbox}

\colorlet{mgu_color}{BrickRed}
\colorlet{gru_color}{RoyalBlue}

\def\boxwidth{0.25}
\def\linewidthlight{0.3 mm}
\def\linewidthheavy{0.6 mm}
\pgfplotsset{
    base boxplot style/.style={
        boxplot,
        boxplot/draw direction=x,
        boxplot/box extend={\boxwidth},
    },
    mgu box style/.style={
        base boxplot style,
        boxplot/every box/.style={
            solid, draw=mgu_color, fill=mgu_color!30, line width={\linewidthlight}
        },
        boxplot/every whisker/.style={
            solid, draw=mgu_color, line width={\linewidthlight}
        },
        boxplot/every median/.style={
            solid, draw=mgu_color, line width={\linewidthheavy}
        },
        every mark/.style={
            solid, mark=x, draw=mgu_color, mark size=2pt, line width={\linewidthlight}
        },
    },
    mgu bar style/.style={
        xbar, 
        bar width=\boxwidth, 
        draw=mgu_color, 
        fill=mgu_color!30, 
        line width=\linewidthlight,
        mark = none,
        nodes near coords, % MOSTRA I NUMERI
        nodes near coords align={horizontal}, % ALLINEA IL TESTO
        point meta=x % Usa il valore X come testo
    },
    gru box style/.style={
        base boxplot style,
        boxplot/every box/.style={
            solid, draw=gru_color, fill=gru_color!30, line width={\linewidthlight}
        },
        boxplot/every whisker/.style={
            solid, draw=gru_color, line width={\linewidthlight}
        },
        boxplot/every median/.style={
            solid, draw=gru_color, line width={\linewidthheavy}
        },
        every mark/.style={
            solid, mark=x, draw=gru_color, mark size=2pt, line width={\linewidthlight}
        },
    },
    gru bar style/.style={
        xbar, 
        bar width=\boxwidth, 
        draw=gru_color, 
        fill=gru_color!30, 
        line width=\linewidthlight,
        mark = none,
        nodes near coords, % MOSTRA I NUMERI
        nodes near coords align={horizontal}, % ALLINEA IL TESTO
        point meta=x % Usa il valore X come testo
    },
    mgu legend/.style={
        draw=mgu_color, 
        fill=mgu_color!30, 
        line width=\linewidthlight
    },
    gru legend/.style={
        draw=gru_color, 
        fill=gru_color!30, 
        line width=\linewidthlight
    },
}

\newtheorem{theorem}{Theorem}
\newtheorem{proposition}{Proposition}
\newtheorem{definition}{Definition}
\newtheorem{assumption}{Assumption}

\newtheorem{remark}{Remark}

\newcommand{\norm}[2]{ \left\Vert {#1} \right\Vert_{#2}}

\newcommand{\infnorm}[1]{ \left\Vert {#1} \right\Vert}
\newcommand{\abs}[1]{ \left| {#1} \right|}

\newcommand{\defeq}{\vcentcolon=} % Definition equality symbol
\newcommand{\quotes}[1]{``{#1}''}
  
\newcommand{\davide}[1]{{\color{green} \textbf{DP:} #1}}
\newcommand{\stefano}[1]{{\color{violet} \textbf{SDC:} #1}}
\newcommand{\mirko}[1]{{\color{red} \hl{\textbf{MM:}} \hl{#1}}}
\newcommand{\tocite}[1]{{\color{red} \textbf{TO CITE}} }

\newcommand{\algorithmhighlight}[1]{{\color{gray} #1}}

\newcommand{\colorrowoftable}{\rowcolor{gray!40}}

\newcommand{\lr}{\operatorname{\textit{\l{}}}}
\newcommand{\hiddenunits}{h}
\newcommand{\layers}{\mathcal{L}}

\newcommand{\sequence}{v}
\newcommand{\totsequence}{V}
\newcommand{\sequencep}{(\sequence)}
\newcommand{\parameters}{\theta}

\newcommand{\MSE}{\textrm{MSE}}
\newcommand{\fit}{\textrm{Fit}}

\newcommand{\loss}[1]{J^{#1}}
\newcommand{\isspenalty}{\rho}
\newcommand{\margin}{\mu}
\newcommand{\projection}[1]{\Pi_{#1}}
\newcommand{\feasibleset}{\mathcal{F}}
\newcommand{\projectionmargin}{\epsilon}

\newcommand{\dISS}{$\delta$ISS}
\newcommand{\realset}{\mathbb{R}}
\newcommand{\naturalset}{\mathbb{N}} % Include lo zero come definito in notation
\newcommand{\posintegerset}{\mathbb{Z}^+} % Solo numeri positivi {1, 2, ...}

\newlength{\customdisplayspace}
\newlength{\customjotspace}
\newenvironment{tequation}
 {\begingroup
 \let\displaystyle\textstyle
 \setlength{\abovedisplayskip}{\customdisplayspace}
 \setlength{\belowdisplayskip}{\customdisplayspace}
 \setlength{\abovedisplayshortskip}{\customdisplayspace}
 \setlength{\belowdisplayshortskip}{\customdisplayspace}
 \equation}
 {\endequation
 \endgroup}
\newenvironment{tequation*}
 {\begingroup
 \let\displaystyle\textstyle
 \setlength{\abovedisplayskip}{\customdisplayspace}
 \setlength{\belowdisplayskip}{\customdisplayspace}
 \setlength{\abovedisplayshortskip}{\customdisplayspace}
 \setlength{\belowdisplayshortskip}{\customdisplayspace}
 \csname equation*\endcsname}
 {\endequation
 \endgroup}

\newenvironment{talign}
 {\begingroup
 \setlength{\abovedisplayskip}{\customdisplayspace}
 \setlength{\belowdisplayskip}{\customdisplayspace}
 \setlength{\abovedisplayshortskip}{\customdisplayspace}
 \setlength{\belowdisplayshortskip}{\customdisplayspace}
 \setlength{\jot}{\customjotspace}
 \align}
 {\endalign
 \endgroup}
\newenvironment{talign*}
 {\begingroup
 \setlength{\abovedisplayskip}{\customdisplayspace}
 \setlength{\belowdisplayskip}{\customdisplayspace}
 \setlength{\abovedisplayshortskip}{\customdisplayspace}
 \setlength{\belowdisplayshortskip}{\customdisplayspace}
  \setlength{\jot}{\customjotspace}
 \csname align*\endcsname}
 {\endalign
 \endgroup}

\newcommand{\maxuandh}{\max_{ { \scriptscriptstyle \boldsymbol{\tilde{u}}_k^{(l)} \in \tilde{\mathcal{U}}^{(l)},\boldsymbol{h}_{k}^{(l)} \in \mathcal{H}_{\mathrm{inv}}^{(l)}} } }
\newcommand{\maxuandhstackedb}{\max_{\scalebox{0.65}{$\begin{array}{c} \boldsymbol{\tilde{u}}_k^{(l), \mathrm{b}} \in \tilde{\mathcal{U}}^{(l)}, \\[-2pt] \boldsymbol{h}_{k}^{(l), \mathrm{b}} \in \mathcal{H}_{\mathrm{inv}}^{(l)}\end{array}$}}}

\allowdisplaybreaks[4]

\begin{document}
\begin{frontmatter}
	\runtitle{Stability properties of Minimal Gated Unit neural networks}
	\title{Stability properties of Minimal Gated Unit Neural Networks\thanksref{footnoteinfo}}

	\thanks[footnoteinfo]{The material in this paper was not presented at any conference.}
	\thanks[corresp]{Corresponding author.}

	\author{Stefano De Carli\thanksref{corresp}}\ead{stefano.decarli@unibg.it},
	\author{Davide Previtali}\ead{davide.previtali@unibg.it},
	\author{Mirko Mazzoleni}\ead{mirko.mazzoleni@unibg.it},
	\author{Fabio Previdi}\ead{fabio.previdi@unibg.it}

	\address{Department of Management, Information and Production Engineering, University of Bergamo, Via G. Marconi 5, 24044 Dalmine, Bergamo, Italy}

	\begin{abstract}
In this work, we address the need for efficient and formally stable Recurrent Neural Networks (RNNs) in environments with limited computational resources by analyzing the stability of the Minimal Gated Unit (MGU) network, a lightweight alternative to common gated RNNs used in system identification.
We derive sufficient parametric conditions 
%, expressed as nonlinear inequalities, 
for the MGU network's input-to-state stability and incremental input-to-state stability properties.
These conditions enable a-posteriori validation of model stability and form the basis for novel stability-promoting training methodologies, including a warm-start of the network's parameters and a projected gradient-based optimization scheme, both of which are presented in this work.
%
% A modified, stability-driven, early stopping strategy is also proposed to account for stability requirements when validating the network's performance during training.
%
Comparative evaluation, including robustness analysis and validation on synthetic and real-world data (i.e., the Silverbox benchmark), demonstrates that the minimal gated unit network successfully combines formal stability guarantees with superior parameter efficiency and faster inference times compared to other state-of-the-art recurrent neural networks, while maintaining comparable and satisfactory accuracy.
Notably, the results attained on the Silverbox benchmark illustrate that the stable MGU network effectively captures the system dynamics, whereas other stable RNNs fail to converge to a reliable model.

%\davide{Metterei frase per valorizzare performance top su Silverbox benchmark}
%
\end{abstract}


\end{frontmatter}

\section{Introduction}
\label{s:introduction}
\emph{Recurrent Neural Networks (RNNs)}, in particular \emph{gated} RNNs~\cite{murphyProbabilisticMachineLearning2025} such as Long Short-Term Memory (LSTM) and Gated Recurrent Unit (GRU) networks, are widely adopted in the control and system identification communities~\cite {ljungDeepLearningSystem2020, pillonettoDeepNetworksSystem2025}.
These networks offer powerful internal memory mechanisms for modeling nonlinear dynamical behaviors in Multiple-Input Multiple-Output (MIMO) systems.
Furthermore, their inherent state-space formulation makes them ideal for integration into model-based control strategies, leading to many recent applications within the Model Predictive Control (MPC) framework~\cite{lawrynczukLSTMGRUType2025,terzi_learning_2021,bonassiRecurrentNeuralNetworks2022,schimpernaRobustConstrainedNonlinear2024, schimpernaRobustOffsetFreeConstrained2024}.

However, despite their proven modeling ability, these gated RNNs are built upon architecturally complex structures~\cite{murphyProbabilisticMachineLearning2025}.
Specifically, LSTM and GRU networks, with their multiple-gate mechanisms, inherently carry a \emph{high parameter count} and \emph{significant computational load} during training and inference.
This complexity poses a major barrier for applications with limited hardware resources, such as those used in embedded control systems~\cite{mayneModelPredictiveControl2014}, where fast inference times and a low memory footprint are important requirements.

Beyond implementation constraints, a significant theoretical challenge lies in the limited formal guarantees available for the identified RNNs~\cite{pillonettoDeepNetworksSystem2025}.
For instance, traditional training of RNNs \emph{does not inherently ensure stability properties}, leading to potential physically-inconsistent behaviors and unbounded predictions~\cite{millerStableRecurrentModels2018a}.
To address this issue, the control community has recently worked on establishing formal stability properties for these networks~\cite{terzi_learning_2021,bonassi_stability_2021,bonassiRecurrentNeuralNetworks2022,decarliInfinitynormbasedInputtoStateStableLong2025}, primarily focusing on \emph{Input-to-State Stability (ISS)} and \emph{Incremental Input-to-State Stability (\dISS)}~\cite{sontagInputtoStateStabilityProperty1995}.
While ISS ensures that bounded inputs lead to bounded state trajectories~\cite{jiang_input--state_2001}, the stricter \dISS{} property guarantees that the distance between any two state trajectories is bounded by the difference between their respective initial states and input sequences~\cite{sontagInputtoStateStabilityProperty1995}.
Establishing these properties is vital for the identification task itself, to have a consistent model when the identified system is itself stable, but it becomes even more critical when the identified models are intended for use in control loops, where one is interested in the best stable model that approximates the given system~\cite[Chapter 12]{borrelliPredictiveControlLinear2017}.
For instance, formal stability is a prerequisite for ensuring the safety of MPC controllers~\cite{borrelliPredictiveControlLinear2017,bonassiRecurrentNeuralNetworks2022} and for designing converging state observers~\cite{schimpernaRobustConstrainedNonlinear2024, schimpernaRobustOffsetFreeConstrained2024}.

Addressing this, \emph{sufficient parametric conditions} have been derived for the \emph{\dISS{} of LSTM~\cite{terzi_learning_2021} and GRU~\cite{bonassi_stability_2021} networks}, paralleling efforts in other model classes such as kernel-based methods~\cite{scandellaKernelBasedIdentificationIncrementally2023,scandellaKernelBasedLearningStable2024} and switched systems~\cite{liuIntegralInputtoStateStabilitySwitched2022}.
This context highlights the need for alternative architectures capable of reconciling the powerful modeling capabilities of gated RNNs with the parameter efficiency required by embedded applications, while simultaneously allowing for the derivation of parametric conditions to promote their stability during training.

%\paragraphfont{Contributions}
%\paragraph*{Contributions}
\subsection{Contributions}
In this work, we introduce the analysis and promotion of stability in \emph{Minimal Gated Unit (MGU) networks}~\cite{zhouMinimalGatedUnit2016}, a \emph{parameter-efficient} and \emph{computationally lighter} alternative to GRU and LSTM networks, along with its \emph{first} validation in system identification tasks.
The primary contributions of this paper are:
\begin{enumerate}
    \item We extend the %infinity-norm 
    stability analysis framework, previously established for GRU networks in~\cite{bonassi_stability_2021}, to the MGU architecture.
    Reflecting the MGU network's specific structure, which induces a tighter coupling between memory retention and input updates compared to the GRU architecture, we derive \emph{distinct parametric conditions that structurally differ} from their GRU network counterparts.
    We establish these sufficient conditions for ISS and \dISS{} in both single-layer and multi-layer configurations, and we prove that the \dISS{} condition parametrically implies the ISS one.
    \item We evaluate three methodologies designed to promote the derived \dISS{} condition during network training.
    While loss augmentation constitutes the state-of-the-art approach currently employed in RNNs to promote stability~\cite{terzi_learning_2021,bonassi_stability_2021,decarliInfinitynormbasedInputtoStateStableLong2025}, we introduce and evaluate two novel strategies: (i) \emph{warm-start} of the parameters via an ad-hoc pre-optimization, and a (ii) \emph{projected gradient-based optimization scheme} utilizing a convex feasible set.
    All methodologies are complemented by a \emph{modified early stopping strategy}, inspired by the ones in~\cite{decarliInfinitynormbasedInputtoStateStableLong2025, bonassi_stability_2021}, that ensures that the best network parameters found during training also imply stability and not just superior performance on a validation dataset.
    \item We perform a numerical evaluation of the MGU network using \emph{three benchmark datasets}: the pH Reactor dataset~\cite{terzi_ph_2020} to compare MGU and GRU networks and to test the proposed training methodologies, the Four-Tank system~\cite{alvaradoComparativeAnalysisDistributed2011a} to assess robustness against measurement noise, and the Silverbox dataset~\cite{wigrenThreeFreeData2013} to validate the performance on real-world data.
    %
    % These analysis can be replicated using the code available at \url{https://github.com/StefanoDeCarli/MGU_dISS.git}.
    %
    \item We provide \emph{open-source} MATLAB code, available at~\cite{githubOfTheMethod}, allowing testing the proposed training methodologies for MGU networks.
\end{enumerate}

%\paragraphfont{Paper Organization}
%\paragraph*{Paper organization}
\subsection{Paper organization}
This paper is organized as follows.
Section~\ref{s:preliminaries} reviews the necessary notation, introduces the problem statement, and provides formal definitions for the considered stability properties. %ISS and \dISS{}. 
Section~\ref{s:mgu} introduces the MGU network architecture. 
Section~\ref{s:stability_mgu} derives the sufficient parametric conditions for the ISS and \dISS{} of both single-layer and multi-layer MGU networks. 
Section \ref{s:stability_promoted_training} details the three methodologies proposed for promoting stability during network training and the stability-driven early stopping strategy. 
Section~\ref{s:numerical_results} provides the numerical evaluation of the MGU network capabilities.
Finally, Section \ref{s:conclusion} summarizes the main findings and results of this work.
\section{Preliminaries and problem statement}
\label{s:preliminaries}
\subsection{Notation and preliminaries}
\label{ss:notation_preliminaries}
Throughout this paper, $\realset$ and $\naturalset$ refer to the sets of real numbers and natural numbers (including zero), respectively.
Positive integers are indicated by $\posintegerset$, %$\posintegerset = \{1, 2, \dots \}$, 
while $\realset_{>0}$ and $\realset_{\geq0}$ stand for the sets of positive and non-negative real numbers, respectively.
Given $n, m \in \posintegerset$, $n \mod m$ is the remainder of the division of $n$ by $m$.
Further, $\realset^n$ is the set of real column vectors of dimension $n$, while $\realset^{n \times m}$ is the set of real matrices of dimension $n \times m$.
The cardinality of a set $\mathcal{S}$ is indicated by $\abs{\mathcal{S}}$.
Vectors are typeset in boldface lowercase; for instance, an $n$-dimensional column vector is written as $\boldsymbol{v}=\left[ v_1,\dots, v_n\right]^\top \in \realset^n$.
Furthermore, $\mathbf{0}_{n} \in \realset^{n}$ and $\mathbf{1}_{n} \in \realset^{n}$ are the $n$-dimensional vectors of zeros and ones, respectively.
%Furthermore, $\boldsymbol{g}_{n} \in \realset^{n}$ refers to a vector with all entries equal to the scalar $g \in \realset$. % is the $n$-dimensional column vector with all entries equal to the scalar $g\in\realset$.
%
Next, the Hadamard (element-wise) product is indicated by $\circ$, whereas $\norm{\cdot}{p}$ represents the $p$-norm of a vector or matrix. 
Unless otherwise specified, $\infnorm{\cdot}$ stands for the infinity norm (i.e., $p=\infty$).

We represent scalar discrete-time signals as $s_k$, where $k \in \naturalset$ is the time step.
Correspondingly, an $n$-dimensional discrete-time signal is defined as $\boldsymbol{x}_k = \left[x_{1, k}, \dots, x_{n, k}\right]^\top \in \realset^n$.
The sequence of values attained by $\boldsymbol{x}_k$ between $k = k_1$ and $k = k_2$ is denoted by %$\mathbf{x} = \{ \boldsymbol{x}_k \}_{k=k_1}^{k_2}$.
$\{ \boldsymbol{x}_k \}_{k=k_1}^{k_2}$.
The sample mean of a finite signal $s_k$ is indicated by $\bar{s}$.

We define the \emph{sigmoid} function as $\sigma (x) \defeq{} \left(1 + e^{-x}\right)^{-1} \in (0,1)$ and the \emph{hyperbolic tangent} function as $\phi(x) \defeq{} \tanh(x) \in (-1,1)$, where $\phi(x)=\left(e^{x} - e^{-x}\right)\left(e^{x} + e^{-x}\right)^{-1}$, both functions $\forall x \in \realset$.
For a scalar function $f:\realset\to\realset$, its bold counterpart $\boldsymbol{f}:\realset^n\to\realset^n$, $n \in \posintegerset$, indicates its component-wise application to a vector, so that, for any $\boldsymbol{v}\in\realset^n,\boldsymbol{f}(\boldsymbol{v})=[f(v_1),\dots,f(v_n)]^\top$.
%
% Regarding neural networks, the superscript $(l)$, e.g. $v^{(l)}$, identifies a quantity referred to the $l$-th layer.

\begin{comment}
Note that $\sigma (x) \in (0,1)$ and $\phi (x) \in (-1,1)$.
%
Both functions are Lipschitz-continuous with respective coefficients $L_\sigma = \frac{1}{4}$ and $L_{\phi} = 1$~\cite{goodfellowDeepLearning2016}.
\end{comment}

We employ the standard definitions of comparison functions~\cite{kellettCompendiumComparisonFunction2014}.
A continuous function $\gamma:\realset_{\geq0}\to\realset_{\geq0}$ belongs to class $\mathcal{K}$ if it is strictly increasing with $\gamma\left(0\right)=0$.
It is a class $\mathcal{K}_{\infty}$ function if, additionally, $\gamma\left(r\right)\to\infty$ as $r\to\infty$.
A continuous function $\psi:\realset_{\geq0}\times\naturalset\to\realset_{\geq0}$ is a class $\mathcal{KL}$ function if, for each fixed $k \in \naturalset$, $\psi(r, k)$ is a class $\mathcal{K}$ function and, for each fixed $r \in \realset_{\geq0}$, $\psi(r, k)$ is strictly decreasing to zero as $k \to \infty$.

\subsection{Problem statement}
\label{ss:problem_statement}
We consider the problem of identifying a nonlinear discrete-time MIMO dynamical system from %experimental 
data.
The objective is to learn the unknown map governing the evolution of the system's output $\boldsymbol{y}_k \in \realset^{n_y}, n_y \in \posintegerset$, from its input $\boldsymbol{u}_k \in \realset^{n_u}, n_u \in \posintegerset$, using a dataset $\mathcal{D} = \{ \mathcal{D}^{\sequencep} \}_{\sequence = 1}^{\totsequence}$ comprising $\totsequence \in \posintegerset$ input-output sequences. 
Each sequence is defined as $\mathcal{D}^{\sequencep} \defeq{} \{(\boldsymbol{u}_k^{\sequencep}, \boldsymbol{y}_k^{\sequencep})\}_{k=0}^{N^{\sequencep}-1}$, $v \in \{1, \ldots, \totsequence\}$, where $N^{\sequencep} \in \posintegerset$ represents the sequence length.
%
% The system is assumed to be characterized by a state $\boldsymbol{x}_k \in \realset^{n_x}, n_x \in \posintegerset$, where the inputs $\boldsymbol{u}_k \in \mathcal{U} \subset \realset^{n_u}$ and the states $\boldsymbol{x}_k \in \mathcal{X} \subset \realset^{n_x}$ belong to the compact sets $\mathcal{U}$ and $\mathcal{X}$, respectively.
%
% The initial state $\boldsymbol{x}_0 \in \mathcal{X}$ is assumed to be unknown and is not subject to estimation.
% However, the state of the system $\boldsymbol{x}_k$ (and its dimension $n_x$) is assumed to be unknown and is not subject to estimation.
%
The system is characterized by a state $\boldsymbol{x}_k \in \realset^{n_x}, n_x \in \posintegerset$, that is assumed to be \emph{unknown} and \emph{unmeasurable}.
Notably, \emph{no stability assumptions are imposed on the underlying system's dynamics}.

Our goal is to \emph{identify the system's input-output dynamics by employing an RNN architecture}, specifically an MGU network.
Explicit identification of the state $\boldsymbol{x}_k$ is \emph{not a concern}.
%; rather, we are only interested in estimating the input-output map.
%
We focus on the classes of MGU models that are input-to-state stable and incremental input-to-state stable.
In the following, the stability properties are defined for a generic system with state equation $\boldsymbol{x}_{k+1} = \boldsymbol{w}(\boldsymbol{x}_k, \boldsymbol{u}_k)$, $\boldsymbol{w}: \realset^{n_x} \times \realset^{n_u} \to \realset^{n_x}$, and stated with respect to the \emph{compact} sets $\mathcal{X} \subset \mathbb{R}^{n_x}$ and $\mathcal{U} \subset \mathbb{R}^{n_u}$, assuming that $\mathcal{X}$ is \emph{forward invariant}, i.e., for any $\boldsymbol{x}_k \in \mathcal{X}$ and $\boldsymbol{u}_k \in \mathcal{U}$, we have $\boldsymbol{w}(\boldsymbol{x}_k, \boldsymbol{u}_k) \in \mathcal{X}$, $\forall k \in \naturalset$.
\begin{definition}[ISS~\cite{jiang_input--state_2001,terzi_learning_2021}]
    \label{def:iss}
    A system $\boldsymbol{x}_{k+1} = \boldsymbol{w}(\boldsymbol{x}_k, \boldsymbol{u}_k)$ is ISS in $\mathcal{X}$ with respect to $\mathcal{U}$ if there exist functions $\psi \in \mathcal{KL}$ and $\gamma_u, \gamma_b \in \mathcal{K}_\infty$ such that, for any $k \in \naturalset$, any initial state $\boldsymbol{x}_0 \in \mathcal{X}$, any input sequence $\{\boldsymbol{u}_z \in \mathcal{U} \}_{z=0}^{k-1}$, and any bias vector $\boldsymbol{b} \in \realset^{n_x}$, it holds that:
    \begin{tequation}
        \label{eq:iss}
        \infnorm{\boldsymbol{x}_k} \leq  \psi(\infnorm{\boldsymbol{x}_0}, k) + \gamma_u \left(\max_{0 \leq z < k} \infnorm{\boldsymbol{u}_z} \right) + \gamma_b \left( \infnorm{\boldsymbol{b}} \right).
    \end{tequation}
\end{definition}
\begin{remark}
    The bias term $\gamma_b \left( \infnorm{\boldsymbol{b}} \right)$ in~\eqref{eq:iss} implies that Definition~\ref{def:iss} formally describes input-to-state practical stability~\cite{mironchenkoCriteriaInputtoStatePractical2019}.
    Standard ISS is recovered when $\boldsymbol{b} = \boldsymbol{0}_{n_x}$, as $\gamma_b(0) = 0$.
    We adopt this practical definition of ISS as the bias term $\gamma_b \left( \infnorm{\boldsymbol{b}} \right)$ is necessary to account for the intrinsic bias vectors present in the employed neural network models.
\end{remark}
\begin{definition}[\dISS{}~\cite{bayerDiscretetimeIncrementalISS2013,terzi_learning_2021}]
    \label{def:diss}
    A system $\boldsymbol{x}_{k+1} = \boldsymbol{w}(\boldsymbol{x}_k, \boldsymbol{u}_k)$ is \dISS{} in $\mathcal{X}$ with respect to $\mathcal{U}$ if there exist functions $\psi_{\delta}\in\mathcal{KL}$ and $\gamma_{\delta u}\in\mathcal{K}_{\infty}$ such that, for any $k\in\naturalset$, any pair of initial states $\boldsymbol{x}_{0}^{\mathrm{a}}, \boldsymbol{x}_{0}^{\mathrm{b}} \in\mathcal{X}$, and any pair of input sequences $\{\boldsymbol{u}_{z}^{\mathrm{a}}\in\mathcal{U}\}_{z=0}^{k-1}$ and $\{\boldsymbol{u}_{z}^{\mathrm{b}}\in\mathcal{U}\}_{z=0}^{k-1}$, the corresponding state trajectories $\boldsymbol{x}_{k}^{\mathrm{a}}$ and $\boldsymbol{x}_{k}^{\mathrm{b}}$ satisfy:
    \begin{tequation}
        \label{eq:diss}
        \infnorm{\boldsymbol{x}_{k}^{\mathrm{a}} \! - \! \boldsymbol{x}_{k}^{\mathrm{b}}} \! \leq \! \psi_{\delta} \! \left(\infnorm{ \boldsymbol{x}_{0}^{\mathrm{a}} \! - \! \boldsymbol{x}_{0}^{\mathrm{b}}}\!, k\right) + \gamma_{\delta u} \! \left(\! \max_{0\leq z<k}\infnorm{ \boldsymbol{u}_{z}^{\mathrm{a}} \! -  \! \boldsymbol{u}_{z}^{\mathrm{b}}} \!\right) \!\!.
    \end{tequation}
\end{definition}
\begin{remark}
    A system that is ISS (respectively, \dISS{}) with respect to one norm is ISS (\dISS{}) with respect to any other norm~\cite{jiang_input--state_2001}, although the associated $\mathcal{KL}$ and $\mathcal{K}_\infty$ functions will differ.
    In this work, we define the ISS and \dISS{} properties using the infinity norm for analytical and practical convenience.
\end{remark}
\section{Minimal gated unit networks}
\label{s:mgu}
An MGU network, introduced in~\cite{zhouMinimalGatedUnit2016} for general machine learning tasks, is a type of recurrent neural network architecture designed to simplify the GRU network~\cite{murphyProbabilisticMachineLearning2025}.
This is achieved by \emph{reducing the number of internal gates}\footnote{Here, we define a gate as any internal sub-structure of an RNN that relies on nonlinear activation functions (e.g., the sigmoid and hyperbolic tangent functions) and is equipped with tunable parameters to govern information flow.}\emph{ from three to two}.
We consider an MGU network composed of $L \in \posintegerset$ stacked \emph{layers}, where each layer $l \in \layers \defeq{} \{1, \ldots, L\}$ contains $n^{(l)}_\hiddenunits \in \posintegerset$ \emph{hidden units}.
In the following, we will use the superscript $(l)$ to refer to quantities related to the $l$-th layer of the network (such as $n^{(l)}_\hiddenunits$).
Specifically, the $l$-th layer of an MGU network represents a state-space discrete-time nonlinear MIMO dynamical system with \emph{hidden state} $\boldsymbol{h}_k^{(l)} \in \realset^{n_{\hiddenunits{}}^{(l)}}$ and layer-specific input $\boldsymbol{\tilde{u}}_k^{(l)} \in \realset^{n_{\tilde{u}}^{(l)}}, n_{\tilde{u}}^{(l)} \in \posintegerset$, defined as: 
\begin{tequation}
    \label{eq:input_for_each_layer}
		\boldsymbol{\tilde{u}}_k^{(l)} \defeq{}
        \begin{cases}
    		\boldsymbol{u}_k & \mathrm{if}\ l = 1,\\
    		\boldsymbol{h}_{k+1}^{(l-1)} & \mathrm{if}\ l \in \layers \setminus \{1\},
		  \end{cases}
\end{tequation}
making $n_{\tilde{u}}^{(1)} = n_u$ and $n_{\tilde{u}}^{(l)} = n_{\hiddenunits}^{(l-1)}$, $\forall l \in \layers \setminus \{1\}$.
The dynamics of the MGU network are governed by a \emph{forget gate} $\boldsymbol{f}_k^{(l)} \in \realset^{n_{\hiddenunits}^{(l)}}$ and a \emph{candidate hidden state} $\boldsymbol{\tilde{h}}_k^{(l)} \in \realset^{n_{\hiddenunits}^{(l)}}$:
\begin{subequations}
    \label{eq:gates_MGU}
    \begin{talign}
        \label{eq:forget_gate}
        \boldsymbol{f}_k^{(l)} &\defeq{} \boldsymbol{\sigma}\left( W_f^{(l)} {\boldsymbol{\tilde{u}}}_k^{(l)}  +  R_f^{(l)} \boldsymbol{h}_k^{(l)}  +  \boldsymbol{b}_f^{(l)} \right),   \\
        \label{eq:candidate_hidden_state}
        \boldsymbol{\tilde{h}}_k^{(l)} &\defeq{} \boldsymbol{\phi}\left( W_{\tilde{h}}^{(l)} {\boldsymbol{\tilde{u}}}_k^{(l)}  +  R_{\tilde{h}}^{(l)} \left( \boldsymbol{f}_k^{(l)}  \circ  \boldsymbol{h}_k^{(l)} \right)  +  \boldsymbol{b}_{\tilde{h}}^{(l)} \right),
    \end{talign}
\end{subequations}
where $W_f^{(l)}, W_{\tilde{h}}^{(l)} \in \realset^{n^{(l)}_{\hiddenunits} \times n_{\tilde{u}}^{(l)}}$ are the \emph{input weights}, $R_f^{(l)}, R_{\tilde{h}}^{(l)} \in \realset^{n^{(l)}_\hiddenunits{} \times n^{(l)}_\hiddenunits{}}$ are the \emph{recurrent weights}, and $\boldsymbol{b}_f^{(l)}, \boldsymbol{b}_{\tilde{h}}^{(l)}  \in \realset^{n^{(l)}_\hiddenunits{}}$ are the \emph{bias vectors}, each for $l \in \layers$.

The state equation of a \emph{single-layer} MGU network is defined as a component-wise convex combination of the hidden state and the candidate hidden state, regulated by the forget gate:
\begin{tequation}
    \label{eq:MGU_layer}
    \boldsymbol{h}_{k+1}^{(l)} \defeq{} \left(\boldsymbol{1}_{n_{\hiddenunits}^{(l)}} - \boldsymbol{f}_{k}^{(l)} \right) \circ \boldsymbol{h}_{k}^{(l)} + \boldsymbol{f}_{k}^{(l)} \circ \boldsymbol{\tilde{h}}_{k}^{(l)}.
\end{tequation}
These single-layer MGU network dynamics are illustrated in \figurename{}~\ref{fig:MGU}.
\begin{figure}[tb]
    \centering
    %--- Dimension setup ------------------------------------
\tikzmath{\lblock = 0.8 cm;
\radius = \lblock;
\eqwidth = 3 cm;
\rnnheight = \lblock;
\grid = 0.4;
\halfgrid = \grid/2;}

%--- Style setup ---------------------------------------
\tikzset{
block/.style={
    rectangle,
    rounded corners,
    draw,
    thick,
    minimum width=\lblock,
    minimum height=\lblock,
    align=center
    },
eqblock/.style={
    rectangle,
    rounded corners,
    draw,
    thick,
    minimum width=\eqwidth,
    minimum height=\rnnheight,
    align=center
    },
surround/.style={
    rectangle,
    draw,
    thick
    },
arrow/.style={
    ->,
    >=stealth,
    thick
    },
line/.style={
    thick
    }
}

%--- The picture ---------------------------------------
\begin{tikzpicture}[every node/.style={font=\small}]
    % Nodes
    % Inputs
    \node[block,draw=blue, text=blue] (uk) {$\boldsymbol{\tilde{u}}_k^{(l)}$};
    \node[block,draw=orange, text=orange, above=\grid cm of uk.north] (hk) {$\boldsymbol{h}_k^{(l)}$};
    % Gates
    \node[block, right=\grid cm of hk.east] (fk) {$\boldsymbol{f}_k^{(l)}$};
    \node[block, right=\grid cm of uk.east] (htk) {$\boldsymbol{\tilde{h}}_k^{(l)}$};
    % Update equation
    \node[eqblock, right=\halfgrid cm - 0.1 cm of htk.east] (eq) at ($ (htk.east)!0.5!(fk.east) $) {$\textcolor{orange}{\boldsymbol{h}_{k+1}^{(l)}} \! = \! (\boldsymbol{1}_{\!n^{(l)}_\hiddenunits} \! - \! \boldsymbol{f}_k^{(l)}) \! \circ \! \textcolor{orange}{\boldsymbol{h}_k^{(l)}} \!\! + \! \boldsymbol{f}_k^{(l)} \!\! \circ \!\boldsymbol{\tilde{h}}_k^{(l)}$};
    % Output
    \node[block,draw=orange, text=orange, right=\halfgrid cm of eq.east] (hkp) {$\boldsymbol{h}_{k+1}^{(l)}$};

    % Connections
    % Input connections
    \draw[arrow,draw=blue] (uk) -- (htk);
    \draw[arrow,draw=blue] ($ (uk.east)!0.5!(htk.west) $) |- (fk.west |- {$(fk.south)!0.25!(fk.north)$});
    % State connections
    \draw[arrow,draw=orange] (hk) -- (fk);
    \draw[arrow,draw=orange] ($ (hk.east)!0.25!(fk.west) $) |- (htk.west |- {$(htk.north)!0.25!(htk.south)$});
    \draw[arrow,draw=orange] ($ (hk.east)!0.25!(fk.west) $) --++ (90:0.65 cm) -| (eq.north -| {$(eq.west)!0.62!(eq.east)$});
    % Gate connections
    \draw[arrow] (fk.south -| {$(fk.west)!0.3!(fk.east)$}) -- (htk.north -| {$(htk.west)!0.3!(htk.east)$});
    \draw[line] (fk.south -| {$(fk.west)!0.7!(fk.east)$}) -- (htk.north -| {$(htk.west)!0.7!(htk.east)$});
    \coordinate (g1) at (fk.south -| {$(fk.west)!0.7!(fk.east)$});
    \coordinate (g2) at (htk.north -| {$(htk.west)!0.7!(htk.east)$});
    \draw[arrow] ($ (g1)!0.5!(g2) $) -- (eq);
    % Output connection
    \draw[arrow] (eq) -- (hkp);

\end{tikzpicture}
    \vspace{-0.2cm}
    \caption{MGU $l$-th layer architecture at time $k$.}
    \label{fig:MGU}
\end{figure}
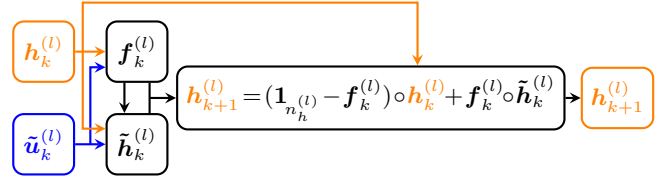
\begin{assumption}[Initial state boundedness]
    \label{as:initial_state_boundedness}
    The initial hidden state of every MGU layer lies within the unit hypercube, i.e., $\boldsymbol{h}_{0}^{(l)} \in \left[-1,1\right]^{n_{\hiddenunits}^{(l)}}$, $\forall l \in \layers$.
\end{assumption}
This initialization is standard practice, often satisfied by setting $\boldsymbol{h}_0^{(l)} = \boldsymbol{0}_{n_{\hiddenunits}^{(l)}}, \forall l \in \layers$~\cite{pascanuHowConstructDeep2013}. 
A key property of the MGU layer is the existence of a compact forward invariant set for the hidden state dynamics, which ensures that the state trajectories remain bounded within the unit hypercube and is essential for the stability properties reviewed in Section~\ref{ss:problem_statement}.
\begin{proposition}[Forward invariant set of the hidden state]
    \label{prop:invariant_state_layer}
    The set $\mathcal{H}_{\mathrm{inv}}^{(l)} \defeq{} \left[-1,1\right]^{n_{\hiddenunits}^{(l)}}$ is a forward invariant compact set for the dynamics in~\eqref{eq:MGU_layer}, meaning that, for any $l \in \layers$, if $\boldsymbol{h}_{0}^{(l)} \in \mathcal{H}_{\mathrm{inv}}^{(l)}$, then $\boldsymbol{h}_{k}^{(l)}\in\mathcal{H}_{\mathrm{inv}}^{(l)}$ for all $k\in\naturalset$ and any $\boldsymbol{\tilde{u}}_k^{(l)} \in \realset^{n_{\tilde{u}}^{(l)}}$.
\end{proposition}
The proof is reported in Appendix~\ref{s:proofs}.

\begin{figure}[tb]
    \centering
    %--- Dimension setup ------------------------------------
\tikzmath{\rnnwidth = 1.2;
\rnnheight = 0.7;
\fcsquare = \rnnheight;
\ioradius = \rnnheight;
\hdist = \rnnwidth/2;
\hdiststate = \rnnwidth*2/3;
\arrowcenterdisp = \rnnwidth/2 - 0.2;}

%--- Style setup ---------------------------------------
\tikzset{
nnbox/.style={
    rectangle,
    rounded corners,
    draw,
    thick,
    minimum width=\rnnwidth cm,
    minimum height=\rnnheight cm,
    align=center,
    fill=lightgray
    },
box/.style={
    rectangle,
    rounded corners,
    draw,
    thick,
    minimum width=\fcsquare cm,
    minimum height=\fcsquare cm,
    align=center,
    fill=lightgray
    },
iobox/.style={
    rectangle,
    rounded corners,
    draw,
    thick,
    minimum width=\fcsquare cm,
    minimum height=\fcsquare cm,
    align=center
    },
arrow/.style={
    ->,
    >=stealth,
    thick
    },
feedbackarrow/.style={
    ->,
    >=stealth,
    thick
    }
}

%--- The picture ---------------------------------------
\begin{tikzpicture}[every node/.style={font=\small}]
    % Nodes
    \node[iobox] (uk) {$\boldsymbol{u}_k$};
    \node[nnbox, right=\hdist/3 of uk] (rnn1) {MGU$^{(1)}$};
    \node[nnbox, right=\hdiststate of rnn1] (rnn2) {MGU$^{(2)}$};
    \node[right= \hdiststate of rnn2] (rnn_empty) {$\cdots$};
    \node[box, right=\hdiststate of rnn_empty] (fc) {FC};
    \node[iobox, right=\hdist/3 of fc] (yk) {$\boldsymbol{y}_k$};

    % Connections
    \draw[arrow] (uk) -- (rnn1);
    \draw[arrow] (rnn1) -- node[above] {$\boldsymbol{h}^{(1)}_{k+1}$} (rnn2);
    \draw[arrow] (rnn2) -- node[above] {$\boldsymbol{h}^{(2)}_{k+1}$} (rnn_empty);
    \draw[arrow] (rnn_empty) -- node[above] {$\boldsymbol{h}^{(L)}_{k+1}$} (fc);
    \draw[arrow] (fc) -- (yk);

    % Feedback loops
    % Feedback for rnn1
    \coordinate (rnn1_start) at ($(rnn1.south) + (-\arrowcenterdisp,0)$);
    \coordinate (rnn1_bend1) at ($(rnn1.south) + (-\arrowcenterdisp,-0.52)$);
    \coordinate (rnn1_bend2) at ($(rnn1.south) + (\arrowcenterdisp,-0.52)$);
    \coordinate (rnn1_end) at ($(rnn1.south) + (\arrowcenterdisp,0)$);
    \draw[feedbackarrow] (rnn1_end) -- (rnn1_bend2) -- node[above=-0.05] {$\boldsymbol{h}^{(1)}_k$} (rnn1_bend1) -- (rnn1_start);
    % Feedback for rnn2
    \coordinate (rnn2_start) at ($(rnn2.south) + (-\arrowcenterdisp,0)$);
    \coordinate (rnn2_bend1) at ($(rnn2.south) + (-\arrowcenterdisp,-0.52)$);
    \coordinate (rnn2_bend2) at ($(rnn2.south) + (\arrowcenterdisp,-0.52)$);
    \coordinate (rnn2_end) at ($(rnn2.south) + (\arrowcenterdisp,0)$);
    \draw[feedbackarrow] (rnn2_end) -- (rnn2_bend2) -- node[above=-0.05] {$\boldsymbol{h}^{(2)}_k$} (rnn2_bend1) -- (rnn2_start);

\end{tikzpicture}
    \vspace{-0.2cm}
    \caption{MGU network with $L$ layers at time $k$.}
    \label{fig:MGU_network}
\end{figure}
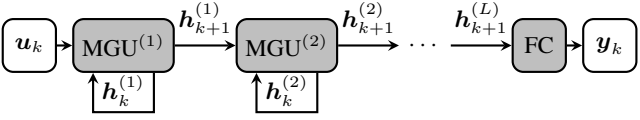

After $L$ stacked MGU layers, a \emph{Fully Connected (FC) layer} generates the network's output $\boldsymbol{y}_k$.
In particular, the overall state-space model of the MGU network reads as:
\begingroup
\setlength{\abovedisplayskip}{\customdisplayspace}
\setlength{\belowdisplayskip}{\customdisplayspace}
\setlength{\abovedisplayshortskip}{\customdisplayspace}
\setlength{\belowdisplayshortskip}{\customdisplayspace}
\setlength{\jot}{\customjotspace}
\begin{subequations}
	\label{eq:MGU_network}
	\begin{numcases}{}
		\label{eq:MGU_network_hidden_state}
		\boldsymbol{h}_{k+1}^{(l)} \! \defeq{} \! \left(\boldsymbol{1}_{n_{\hiddenunits}^{(l)}} \! - \! \boldsymbol{f}_{k}^{(l)}\right) \! \circ \! \boldsymbol{h}_{k}^{(l)} \! + \! \boldsymbol{f}_{k}^{(l)} \! \circ \! \boldsymbol{\tilde{h}}_{k}^{(l)}, \quad \forall l \in \layers, \\
		\label{eq:MGU_network_output}
		\boldsymbol{y}_k \defeq{} W_y \boldsymbol{h}_{k+1}^{(L)} + \boldsymbol{b}_y, 
	\end{numcases}
\end{subequations}
\endgroup
where~\eqref{eq:MGU_network_output} describes the FC layer %mapping the final hidden state $\boldsymbol{h}_{k+1}^{(L)}$ to the output vector $\boldsymbol{y}_k$, 
with $W_y \in \realset^{n_y \times n_{\hiddenunits}^{(L)}}$ and $\boldsymbol{b}_y \in \realset^{n_y}$ being the \emph{output weight matrix} and \emph{output bias vector}, respectively.
The architecture of the multi-layer MGU network, resulting from the stacking of $L$ layers following the interconnection logic in \eqref{eq:input_for_each_layer}, is depicted in \figurename{}~\ref{fig:MGU_network}.

From Proposition~\ref{prop:invariant_state_layer}, it follows that the overall MGU network state $\boldsymbol{h}_k \defeq{} \left[\boldsymbol{h}_{k}^{(1)^\top}, \dots, \boldsymbol{h}_{k}^{(L)^\top}\right]^\top$ in~\eqref{eq:MGU_network_hidden_state} evolves within the Cartesian product $\mathcal{H}_{\mathrm{inv}} \defeq{} \bigtimes_{l=1}^L \mathcal{H}_{\mathrm{inv}}^{(l)} \subset \realset^{n_\hiddenunits{}}$, where $n_\hiddenunits{} \defeq{} \sum_{l=1}^L n_\hiddenunits^{(l)}$ is the total number of hidden units in the network.
Clearly, $\mathcal{H}_{\mathrm{inv}}$ is also \emph{compact} and \emph{forward invariant} as a consequence of Proposition~\ref{prop:invariant_state_layer}.

Overall, the MGU network in~\eqref{eq:MGU_network} relies on the following set of \emph{tunable parameters} (see~\eqref{eq:gates_MGU} and~\eqref{eq:MGU_network_output}):
\begin{subequations}
    \label{eq:parameters}
    \begin{talign}
        \label{eq:MGU_layer_parameters}
        \parameters^{(l)} &\defeq{} \{W_f^{(l)}, W_{\tilde{h}}^{(l)}, R_f^{(l)}, R_{\tilde{h}}^{(l)}, \boldsymbol{b}_f^{(l)}, \boldsymbol{b}_{\tilde{h}}^{(l)}\}, \quad l \in \layers \\
        \label{eq:MGU_parameters}
        \parameters &\defeq{} \{W_y, \boldsymbol{b}_y\}  \cup \bigcup_{l \in \layers} \parameters^{(l)}.
    \end{talign}
\end{subequations}
Then, the total number of parameters is given by:
\begin{tequation}
    \label{eq:MGU_cardinality}
    n_{\theta} = n_y (n_{\hiddenunits}^{(L)} + 1) + 2 \sum_{l = 1}^{L} n^{(l)}_{\hiddenunits}\left( n_{\tilde{u}}^{(l)} + n^{(l)}_{\hiddenunits} + 1 \right).
\end{tequation}
\begin{remark}
    \label{rem:mgu_smaller}
    Given a hidden state dimensionality $n^{(l)}_\hiddenunits{}$, an MGU layer requires two-thirds the parameters of a GRU layer, owing to its two-gate architecture versus the GRU's three-gate architecture\footnote{The number of parameters for a GRU network uses a coefficient of $3$ instead of $2$ for the second addend in~\eqref{eq:MGU_cardinality}.}.
\end{remark}
Having fully characterized the MGU network architecture, 
% including its state representation, dynamics, and intrinsic invariance properties, 
we now proceed to establish its formal stability properties.
\section{Stability properties of MGU networks}
\label{s:stability_mgu}
In this Section, we derive formal stability properties for MGU networks.
First, in Section~\ref{ss:stability_mgu_layers}, we establish \emph{sufficient parametric conditions} for the stability of MGU layers; subsequently, in Section~\ref{ss:stability_mgu_networks}, we extend these findings to ensure the stability of the overall MGU network. 
To carry out this analysis, we make the following Assumption about the inputs of the network.
\begin{assumption}[Input boundedness]
    \label{as:input_boundedness}
    The inputs $\boldsymbol{u}_{k}$ of the MGU network in~\eqref{eq:MGU_network} are bounded such that $\infnorm{\boldsymbol{u}_k} \leq 1, \forall k\in\naturalset$, i.e., $\boldsymbol{u}_k \in \mathcal{U} \defeq{} [-1, 1]^{n_u}$.
\end{assumption}
Assumption~\ref{as:input_boundedness} is standard in neural network applications~\cite{goodfellowDeepLearning2016} and can be readily satisfied by \emph{normalizing} the input vector.
We also remark that similar assumptions were made in prior works on the stability of gated recurrent neural networks, see, e.g., \cite[Assumption 1]{terzi_learning_2021} and~\cite[Assumption 1]{bonassi_stability_2021}.

Concerning the inputs of the $l$-th MGU layer in~\eqref{eq:input_for_each_layer}, $l \in \layers \setminus \{1\}$, note that under Assumption~\ref{as:initial_state_boundedness} and due to Proposition~\ref{prop:invariant_state_layer} we have $\infnorm{\boldsymbol{\tilde{u}}_k^{(l)}} \leq 1$, $\forall k \in \naturalset$.
Then, considering all the inputs in~\eqref{eq:input_for_each_layer}, we can define their respective input spaces as $\tilde{\mathcal{U}}^{(l)} \defeq{} [-1, 1]^{n_{\tilde{u}}^{(l)}}$, $\forall l \in \layers$, implying that $\boldsymbol{\tilde{u}}_k^{(l)} \in \tilde{\mathcal{U}}^{(l)}$, $\forall k \in \naturalset$.
Clearly, all the $\tilde{\mathcal{U}}^{(l)}$'s are \emph{compact} sets, consistent with Section~\ref{ss:problem_statement}.
Further, $\tilde{\mathcal{U}}^{(1)} = \mathcal{U}$ as in Assumption~\ref{as:input_boundedness}.

\subsection{Stability of MGU layers}
\label{ss:stability_mgu_layers}
Given the previous assumptions, we proceed to derive sufficient parametric conditions for ISS and \dISS{} of MGU layers.
\begin{theorem}[ISS of MGU layers]
    \label{th:MGU_layer_ISS}
    Consider the $l$-th, $l \in \layers$, MGU layer in~\eqref{eq:MGU_layer} with parameters $\theta^{(l)}$ in~\eqref{eq:MGU_layer_parameters}.
    Under Assumption~\ref{as:initial_state_boundedness} and Assumption~\ref{as:input_boundedness}, the layer is ISS in $\mathcal{H}_{\mathrm{inv}}^{(l)}$ with respect to $\tilde{\mathcal{U}}^{(l)}$ as in Definition~\ref{def:iss} if:
    \begin{tequation}
        \label{eq:MGU_layer_ISS_condition}
        \bar{\sigma}_{f}^{(l)}\infnorm{ R_{\tilde{h}}^{(l)}}<1,
    \end{tequation}
    where $\bar{\sigma}_{f}^{(l)}$ is the upper bound of the forget gate $\boldsymbol{f}_k^{(l)}$:
    \begin{tequation*}
        \bar{\sigma}_{f}^{(l)} \defeq \sigma\left(\infnorm{ \begin{bmatrix}W_{f}^{(l)} & R_{f}^{(l)} & \boldsymbol{b}_{f}^{(l)}\end{bmatrix}}\right). \\
    \end{tequation*}
\end{theorem}
The proof is reported in Appendix~\ref{s:proofs}.
The corresponding $\mathcal{KL}$ and $\mathcal{K}_{\infty}$ functions are shown in~\eqref{eq:KL_ISS_layer-wise}~-~\eqref{eq:Kinfty_ISS_layer-wise_bias}.
\begin{theorem}[\dISS{} of MGU layers]
    \label{th:MGU_layer_dISS}
    Consider the $l$-th, $l \in \layers$, MGU layer in~\eqref{eq:MGU_layer} with parameters $\theta^{(l)}$ in~\eqref{eq:MGU_layer_parameters}.
    Under Assumption~\ref{as:initial_state_boundedness} and Assumption~\ref{as:input_boundedness}, the layer is \dISS{} in $\mathcal{H}_{\mathrm{inv}}^{(l)}$ with respect to $\tilde{\mathcal{U}}^{(l)}$ as in Definition~\ref{def:diss} if:
    \begin{tequation}
        \label{eq:MGU_layer_dISS_condition}
        \bar{\sigma}_{f}^{(l)} \! + \! \bar{\sigma}_{f}^{(l)^{2}} \! \infnorm{ R_{\tilde{h}}^{(l)}} \! + \! \frac{1}{4} \infnorm{R_{f}^{(l)}} \! \left(\bar{\sigma}_{f}^{(l)} \infnorm{R_{\tilde{h}}^{(l)}} \! + \! \bar\phi_{\tilde{h}}^{(l)} \! + \! 1 \! \right) \!  <  \! 1,
    \end{tequation}
    where $\bar{\sigma}_{f}^{(l)}$ is the forget gate $\boldsymbol{f}_k^{(l)}$ upper bound from Theorem~\ref{th:MGU_layer_ISS} and $\bar\phi_{\tilde{h}}^{(l)}$ is the upper bound of the candidate hidden state $\boldsymbol{\tilde{h}}_k^{(l)}$:
    \begin{tequation*}
        \bar\phi_{\tilde{h}}^{(l)} \defeq  \phi\left(\infnorm{ \begin{bmatrix}W_{\tilde{h}}^{(l)} & R_{\tilde{h}}^{(l)} & \boldsymbol{b}_{\tilde{h}}^{(l)}\end{bmatrix}} \right). \\
    \end{tequation*}
\end{theorem}
The proof is reported in Appendix~\ref{s:proofs}.
The corresponding $\mathcal{KL}$ and $\mathcal{K}_{\infty}$ functions are shown in~\eqref{eq:KL_dISS_layer-wise}~and~\eqref{eq:Kinfty_dISS_layer-wise_input}.
\begin{proposition}
    \label{pr:MGU_layer_dISS_imply_ISS}
    If the $l$-th, $l \in \layers$, MGU layer in~\eqref{eq:MGU_layer} satisfies the sufficient \dISS{} condition in~\eqref{eq:MGU_layer_dISS_condition}, it also satisfies the ISS condition in~\eqref{eq:MGU_layer_ISS_condition}.
\end{proposition}
The proof is reported in Appendix~\ref{s:proofs}.
%
\begin{comment}
\begin{remark}[Structural Trade-off]
    \label{rem:structural_tradeoff}
    It is worth noting that the MGU's efficiency comes with a structural trade-off regarding stability. 
    %
    In a GRU, the multiple gates provide independent degrees of freedom to satisfy the stability condition~\cite{bonassi_stability_2021}. 
    %
    In contrast, the MGU couples these mechanisms into less gates. 
    %
    Consequently, conditions~\eqref{eq:MGU_layer_ISS_condition} and~\eqref{eq:MGU_layer_dISS_condition} impose tighter constraints on the feasible parameter space compared to the GRU~\cite{bonassi_stability_2021}. 
    %
    This reduced flexibility makes finding a stable solution via standard unconstrained gradient descent more challenging, motivating the specific stability-promoting training strategies detailed in Section~\ref{s:stability_promoted_training}.
\end{remark}
\end{comment}
%
\subsection{Stability of MGU networks}
\label{ss:stability_mgu_networks}
We now generalize the layer-level stability conditions established in Section~\ref{ss:stability_mgu_layers} to multi-layer MGU networks.
\begin{theorem}[ISS of MGU networks]
    \label{th:MGU_network_ISS}
    Consider the MGU network governed by the dynamics in~\eqref{eq:MGU_network} and parameters $\theta$ in~\eqref{eq:MGU_parameters}, composed of $L$ stacked MGU layers.
    Under Assumption~\ref{as:initial_state_boundedness} and Assumption~\ref{as:input_boundedness}, the MGU network is ISS in $\mathcal{H}_{\mathrm{inv}}$ with respect to $\mathcal{U}$ as in Definition~\ref{def:iss} if the single-layer ISS condition in~\eqref{eq:MGU_layer_ISS_condition} holds for every layer $l \in \layers$.
\end{theorem}
The proof is reported in Appendix~\ref{s:proofs}.
\begin{theorem}[\dISS{} of MGU networks]
    \label{th:MGU_network_dISS}
     Consider the MGU network governed by the dynamics in~\eqref{eq:MGU_network} and parameters $\theta$ in~\eqref{eq:MGU_parameters}, composed of $L$ stacked MGU layers.
     Under Assumption~\ref{as:initial_state_boundedness} and Assumption~\ref{as:input_boundedness}, the MGU network is \dISS{} in $\mathcal{H}_{\mathrm{inv}}$ with respect to $\mathcal{U}$ as in Definition~\ref{def:diss} if the single-layer \dISS{} condition in~\eqref{eq:MGU_layer_dISS_condition} holds for every layer $l \in \layers$.
\end{theorem}
The proof is reported in Appendix~\ref{s:proofs}.
The corresponding $\mathcal{KL}$ and $\mathcal{K}_{\infty}$ functions are shown in~\eqref{eq:KL_dISS_network-wise}~and~\eqref{eq:Kinfty_dISS_network-wise}.
\begin{proposition}
    \label{pr:MGU_network_dISS_imply_ISS}
    If the MGU network in~\eqref{eq:MGU_network} is \dISS{} as outlined in Theorem~\ref{th:MGU_network_dISS}, then it also is ISS.
\end{proposition}
The proof is reported in Appendix~\ref{s:proofs}. 
%\stefano{Sottolineare che le condizioni sono ottenute considerando i massimi possibili input e stati nei loro set compatti.}
%
%\stefano{Remark per menzionare che queste condizioni scalano male, come quelle per LSTM e GRU in~\cite{terzi_learning_2021,bonassi_stability_2021}, con l'aumentare del numero di hidden unit?}

Theorem~\ref{th:MGU_network_ISS} and Theorem~\ref{th:MGU_network_dISS} establish the sufficient parametric conditions that formally guarantee the stability of the MGU network.
We now propose several methodologies that favor the satisfaction of these stability constraints during the model identification process.
\section{Stability-promoted training of MGU networks}
\label{s:stability_promoted_training}
Estimation of the parameters $\theta$ of the MGU network in~\eqref{eq:MGU_parameters} (more commonly referred to as \emph{MGU network training}) involves solving an \emph{optimization problem} using \emph{gradient-based} optimization procedures specifically tailored for neural networks~\cite{pillonettoDeepNetworksSystem2025}.
Traditional training methodologies (reviewed in Section~\ref{ss:training}) entail the \emph{unconstrained} minimization of a \emph{loss function} that measures the discrepancy between the output predictions of the RNN and the available output data, possibly with the addition of a regularization term~\cite{pillonettoDeepNetworksSystem2025}.
However, this rationale offers \emph{no stability guarantees} for the identified model.
To address this, we introduce and compare three \emph{complementary} methodologies for training stable MGU networks.
These include: (i) the current state-of-the-art approach known as \emph{loss augmentation}~\cite{terzi_learning_2021,bonassi_stability_2021,decarliInfinitynormbasedInputtoStateStableLong2025} (detailed in Section~\ref{ss:loss_augmentation}), which we find insufficient alone for MGU network training; (ii) a parameters \emph{warm-start procedure} (Section~\ref{ss:parameters_warm_start}); and (iii) a \emph{Projected Gradient-based optimization Method (PGM)} (Section~\ref{ss:projected_gradient}).
All methodologies are enriched with a \emph{stability-driven early stopping} strategy, inspired by those in~\cite{decarliInfinitynormbasedInputtoStateStableLong2025, bonassi_stability_2021}, which ensures that MGU network training returns the best-performing \dISS{} network, if available~\footnote{
In our empirical results in Section~\ref{s:numerical_results}, when employing any of the proposed stability-promoting methodologies, we were always able to find a \dISS{} network. 
Nonetheless, if that were not the case, we suggest trying a different network initialization, increasing the stability penalty weight presented in Section~\ref{ss:loss_augmentation}, or changing the MGU network hyperparameters in Algorithm~\ref{alg:MGU_training}.
} (see Section~\ref{ss:loss_augmentation}).

We remark that the proposed methodologies are not mutually exclusive alternatives but can be used in \emph{synergy}, each promoting stability in a different way. 
For clarity, Algorithm~\ref{alg:MGU_training} summarizes the training process of \dISS{} MGU networks, highlighting how the proposed approaches (in gray) interact with the traditional training methodology.
The user may leverage \emph{any or all} the techniques in Sections~\ref{ss:loss_augmentation}-\ref{ss:projected_gradient} based on the application.
% 
%Please refer to Section~\ref{ss:stability_comparison} for further details on how our strategies combine.

\subsection{Traditional recurrent neural network training}
\label{ss:training}
\begin{algorithm*}[tb]
	\caption{MGU network training. 
	\\
	The Algorithm reports the proposed stability-promoting strategies with a \algorithmhighlight{gray} font, which can be enabled or disabled as desired.
	}
	\label{alg:MGU_training}
	\textbf{Inputs}: Training data~$\mathcal{D}_{\mathrm{tr}}$, validation data~$\mathcal{D}_{\mathrm{val}}$, number of mini-batches~$N_{\mathrm{mb}}$, number of layers~$L$ of the MGU network and number of hidden units~$n_{\hiddenunits{}}^{(l)}$ for every layer $ l \in L$, learning rate~$\lr$, learning rate decay factor~$\zeta$, learning rate decay frequency~$e_{\mathrm{dc}}$, maximum epochs~$e_{\mathrm{max}}$, validation frequency~$\kappa_{\mathrm{val}}$, dropout rate~$\xi$, \algorithmhighlight{stability penalty weight $\isspenalty$}, \algorithmhighlight{stability safety margin $\margin$}, \algorithmhighlight{projection margin $\projectionmargin$}.
	\\
	\textbf{Output}: Best model parameters~$\theta^\star$.
	\hrule
	\begin{algorithmic}[1]
		\State Randomly divide the training data $\mathcal{D}_{\mathrm{tr}}$ into $N_{\mathrm{mb}}$ mini-batches $\mathcal{D}_{\mathrm{tr}}^{(b)}$, $b \in \{1, \ldots, N_{\mathrm{mb}}\}$, of roughly the same size
		\State Initialize the MGU network parameters in~\eqref{eq:parameters} via the standard initialization scheme (Section~\ref{ss:training}) to get $\theta_0$
		\State \algorithmhighlight{Project $\theta_0$ onto a stable \dISS{}-compliant region by solving Problem~\eqref{eq:warm_start_optimization}, $\forall l \in \layers$, to get $\tilde{\theta}_0$. Then, replace $\theta_0 \gets \tilde{\theta}_0$ (\emph{parameters warm-start})}
    \State Set $\theta^\star \gets \theta_0$ and ${\loss{}}^\star_\mathrm{val} \gets \infty$
		\State Initialize the epoch and iteration counters: $e \gets 0$ and $\kappa \gets 0$
		\While{$e < e_{\mathrm{max}}$}
			\State Increase the epoch counter $e \gets e + 1$
			\State Group the mini-batches to be used, i.e. $\mathcal{B} \gets \{\mathcal{D}_{\mathrm{tr}}^{(b)}: b \in \{1, \ldots, N_{\mathrm{mb}}\}\}$, randomizing their order
			\If{$e \mod e_{\mathrm{dc}} = 0$}
				\State Apply learning rate decay: $\lr \gets \zeta \lr$
			\EndIf
			\For{$b = 1$ to $N_{\mathrm{mb}}$}
				\State Increase the iteration counter $\kappa \gets \kappa + 1$
				\State Extract the mini-batch $\mathcal{D}_{\mathrm{tr}}^{(b)}$ from $\mathcal{B}$
        \State Enable dropout~\cite{srivastava2014dropout} with probability $\xi$
    		\State Compute the loss function $\loss{}(\theta_{\kappa - 1}; \mathcal{D}_{\mathrm{tr}}^{(b)})$ either as in~\eqref{eq:mse} (traditional training) \algorithmhighlight{or as in~\eqref{eq:augmented_loss} (\emph{loss augmentation})}
    		\State Compute the gradient $\nabla \loss{}(\theta_{\kappa - 1}; \mathcal{D}_{\mathrm{tr}}^{(b)})$
    		\State Update the parameters $\theta_{\kappa - 1}$ as in~\eqref{eq:NN_optimizationStep} to get $\theta_{\kappa}$
				%\State \algorithmhighlight{Project the obtained parameters $\theta_{\kappa}$ onto the feasible set $\feasibleset$ in~\eqref{eq:iss_feasible_set}, i.e., $\theta_{\kappa} \gets \projection{{\feasibleset}} \left(\theta_{\kappa}\right)$ (\emph{projected gradient-based optimization method})}
				\State \algorithmhighlight{Project the parameters onto the feasible set in~\eqref{eq:iss_feasible_set}: $\theta_{\kappa} \gets \projection{{\feasibleset}} \left(\theta_{\kappa}\right)$ (\emph{projected gradient-based optimization method})}
    		\If{$\kappa \mod\kappa_\mathrm{val} = 0$}
            \State Disable dropout~\cite{srivastava2014dropout}
        		\State Compute ${\loss{}}_\mathrm{val}(\theta_{\kappa}; \mathcal{D}_{\mathrm{val}}) \gets \mathrm{\MSE{}}(\theta_{\kappa}; \mathcal{D}_{\mathrm{val}})$ as in~\eqref{eq:mse}
        		\If{${\loss{}}_\mathrm{val}(\theta_{\kappa}; \mathcal{D}_{\mathrm{val}}) < {\loss{}}^\star_\mathrm{val}$ \algorithmhighlight{and Condition~\eqref{eq:MGU_layer_dISS_condition} holds for every $\theta_{\kappa}^{(l)}$, $l \in \layers$ (\emph{stability-driven early stopping})}}
            		\State Set ${\loss{}}^\star_\mathrm{val} \gets {\loss{}}_\mathrm{val}(\theta_{\kappa}; \mathcal{D}_{\mathrm{val}})$
            		\State Update $\theta^\star \gets \theta_{\kappa}$
    		    \EndIf
        \EndIf
    	\EndFor
		\EndWhile
		\State \Return $\theta^\star$
	\end{algorithmic}
\end{algorithm*}
As customary, the data ($\mathcal{D}$) introduced in Section~\ref{ss:problem_statement} is partitioned into three mutually exclusive training ($\mathcal{D}_{\mathrm{tr}}$), validation ($\mathcal{D}_{\mathrm{val}}$), and test ($\mathcal{D}_{\mathrm{tst}}$) datasets, serving for MGU parameters estimation, hyperparameters tuning and early stopping~\cite[Chapter 8]{goodfellowDeepLearning2016}, and final performance evaluation, respectively.
The number of input-output sequences in each dataset is denoted by $V_{\mathrm{tr}},V_{\mathrm{val}},V_{\mathrm{tst}} \in \posintegerset$, respectively, which are such that $V_{\mathrm{tr}} + V_{\mathrm{val}} + V_{\mathrm{tst}} = V$.
%
%Each $v$-th sequence retains its temporal length $N^{(v)}$, as defined in Section 2.2.
%
To ease the computational burden during neural network training, $\mathcal{D}_{\mathrm{tr}}$ is further divided (randomly) into $N_{\mathrm{mb}} \in \posintegerset$ disjoint subsets (the so-called \emph{mini-batches}) $\mathcal{D}_{\mathrm{tr}}^{(b)}$, $b \in \{1, \ldots, N_{\mathrm{mb}}\}$, each containing roughly the same number of input-output sequences and such that $\sum_{b = 1}^{N_{\mathrm{mb}}} \left|\mathcal{D}_{\mathrm{tr}}^{(b)}\right| = V_{\mathrm{tr}}$~\cite[Chapter 8]{goodfellowDeepLearning2016}.

Traditional training of MGU networks (or RNNs in general) involves the minimization of the \emph{Mean Squared Error (MSE) loss function} computed over the training dataset~\cite{pillonettoDeepNetworksSystem2025}:
\begin{talign}
    \label{eq:mse}
    &\loss{}(\parameters; \mathcal{D}_\mathrm{tr}) = \MSE{}(\parameters; \mathcal{D}_{\mathrm{tr}}) \nonumber \\
     & \quad \defeq{} \! \frac{1}{|\mathcal{D}_{\mathrm{tr}}|} \!\!\sum_{\mathcal{D}^{\sequencep} \in \mathcal{D}_{\mathrm{tr}}} \!\!\left[ \!\frac{1}{N^{\sequencep}} \!\!\!\!\! \sum_{k=0}^{N^{\sequencep}-1}\norm{\boldsymbol{y}_k^{\sequencep}-\boldsymbol{\hat{y}}_k^{\sequencep}(\parameters) }{2}^2 \!\right] \!, 
\end{talign}
where $\boldsymbol{\hat{y}}_k^{\sequencep}(\parameters)$ denotes the output predicted by the MGU network in~\eqref{eq:MGU_network} with parameters $\theta$ as in~\eqref{eq:MGU_parameters} for the $v$-th sequence in $\mathcal{D}_\mathrm{tr}$ with inputs $\boldsymbol{u}_k^{\sequencep}$ at the time $k$.
Typically, the minimization of~\eqref{eq:mse} is done using \emph{iterative} \emph{gradient-based} optimizers~\cite{pillonettoDeepNetworksSystem2025}~\cite[Chapter 8]{goodfellowDeepLearning2016}.
Thus, the MGU network parameters in~\eqref{eq:MGU_parameters} are initialized using a \emph{standard initialization scheme} that promotes successful information propagation and breaks symmetry between hidden units, effectively aiding the optimization procedure~\cite[Chapter 8.4]{goodfellowDeepLearning2016}.
In particular, for each layer $l \in \layers$, the MGU network parameters are generally initialized as follows. 
Input and output weight matrices $W_f^{(l)}, W_{\tilde{h}}^{(l)}, W_y$ use Glorot uniform initialization, recurrent weight matrices $R_f^{(l)}, R_{\tilde{h}}^{(l)}$ use orthogonal initialization, and the bias vectors $\boldsymbol{b}_f^{(l)}, \boldsymbol{b}_y$ are typically zero vectors.
The only exceptions are the forget gate bias vectors $\boldsymbol{b}_f^{(l)}$, which are standardly initialized to $\boldsymbol{1}_{n^{(l)}_\hiddenunits{}}$. 
In the following, we mark the parameters in~\eqref{eq:parameters} at the training iteration $\kappa \in \naturalset$ as $\parameters_{\kappa}$ (and $\parameters_\kappa^{(l)}$, $l \in \layers$, for the layer-specific parameters).
Then, the initialization scheme generates $\parameters_0$ (including $\parameters_0^{(l)}$, $\forall l \in \layers$).

At each successive iteration ($\kappa \geq 1$), a gradient-based optimizer updates the MGU network parameters $\parameters_{\kappa - 1}$ using an \emph{update rule} that leverages the gradient $\nabla \loss{}(\theta_{\kappa - 1}; \mathcal{D}_{\mathrm{tr}}^{(b)})$ of the loss function in~\eqref{eq:mse} computed using only a mini-batch of data $\mathcal{D}_{\mathrm{tr}}^{(b)}$, $b \in \{1, \ldots, N_{\mathrm{mb}}\}$, and based on a \emph{learning rate} $\lr \in \mathbb{R}_{> 0}$ that regulates the magnitude of the update~\cite[Chapter 8]{goodfellowDeepLearning2016}:\footnote{We omit other possible hyperparameters (see~\cite[Chapter 8]{goodfellowDeepLearning2016}) since the learning rate is usually the most impactful.}
\begin{tequation}
    \label{eq:NN_optimizationStep}
	\theta_{\kappa}
	=
	\theta_{\kappa - 1} - \mathrm{GradStep}\left(\nabla \loss{}(\theta_{\kappa - 1}; \mathcal{D}_{\mathrm{tr}}^{(b)}), \lr \right),
\end{tequation}
where $\mathrm{GradStep}(\nabla \loss{}(\theta_{\kappa - 1}; \mathcal{D}_{\mathrm{tr}}^{(b)}), \lr )$ is a set of operations that depends on the employed optimization method.
For example, in the case of stochastic gradient descent~\cite[Chapter 8]{goodfellowDeepLearning2016}, we simply have $\mathrm{GradStep}(\nabla \loss{}(\theta_{\kappa - 1}; \mathcal{D}_{\mathrm{tr}}^{(b)}), \lr ) \! = \! \lr \nabla \loss{}(\theta_{\kappa - 1}; \mathcal{D}_{\mathrm{tr}}^{(b)})$. 
Nonetheless, other, more advanced, optimization strategies such as Adam and RMSProp exist~\cite[Chapter 8]{goodfellowDeepLearning2016}.
The mini-batches for the update rule in~\eqref{eq:NN_optimizationStep} are processed sequentially, so that traversing the entire training dataset $\mathcal{D}_{\mathrm{tr}}$ requires $N_{\mathrm{mb}}$ iterations.
When that happens, we say that a training \emph{epoch} has passed.

The optimization procedure is stopped once a maximum number of epochs $e_{\mathrm{max}} \in \posintegerset$ is reached.
The learning rate is decayed by a factor $\zeta \in (0,1)$ every $e_{\mathrm{dc}} \in \posintegerset$ epochs to aid parameters convergence.
To prevent overfitting and improve generalization, two implicit regularization strategies~\cite{pillonettoDeepNetworksSystem2025} are employed: \emph{dropout}~\cite{srivastava2014dropout} with a probability $\xi \in [0,1)$ and \emph{early stopping}~\cite[Chapter 8]{goodfellowDeepLearning2016}.
Concerning early stopping, we compute the mean squared error in~\eqref{eq:mse} on the validation dataset $\mathcal{D}_{\mathrm{val}}$, i.e. $\MSE{}(\parameters_{\kappa}; \mathcal{D}_{\mathrm{val}})$, every $\kappa_{\mathrm{val}} \in \posintegerset$ iterations, and the parameter set $\parameters^\star$ yielding the best validation performance is selected as the optimal one.
%
%Algorithm~\ref{alg:MGU_training} summarizes the training procedure for MGU networks.

\subsection{Loss augmentation and stability-driven early stopping}
\label{ss:loss_augmentation}
The traditional RNN training methodology described in Section~\ref{ss:training} \emph{does not guarantee} that the MGU layers will intrinsically satisfy the ISS condition in~\eqref{eq:MGU_layer_ISS_condition} or the $\delta$ISS condition in~\eqref{eq:MGU_layer_dISS_condition}.
To address this, we augment the loss function in~\eqref{eq:mse} with a penalty term that discourages the violation of these stability conditions.
This approach is inspired by previous works on RNN stability~\cite{terzi_learning_2021, bonassi_stability_2021, decarliInfinitynormbasedInputtoStateStableLong2025}.
Since Proposition~\ref{pr:MGU_layer_dISS_imply_ISS} and Proposition~\ref{pr:MGU_network_dISS_imply_ISS} demonstrate that the \dISS{} condition implies the ISS one, both on the layer and network level, we focus our strategy exclusively on satisfying the stricter \dISS{} requirement.
Then, the proposed \emph{augmented loss function} reads as:
\begin{talign}
    \label{eq:augmented_loss}
        \loss{}(\parameters; \mathcal{D}_\mathrm{tr}) &\defeq \MSE{}(\parameters; \mathcal{D}_\mathrm{tr}) +  \\
        &\quad + \isspenalty \sum_{l \in \layers} \max\left\{ \delta\mathrm{ISS}^{(l)}(\parameters^{(l)}) - (1 - \margin), 0  \right\}, \nonumber
\end{talign}
where $\isspenalty \in \realset_{\geq 0}$ is the \emph{penalty weight} and $\margin \in (0, 1)$ is a \emph{safety margin}, which are tuned empirically.
The \emph{penalty term} $\delta\mathrm{ISS}^{(l)}(\parameters^{(l)})\in \realset_{\geq 0}$ is derived from~\eqref{eq:MGU_layer_dISS_condition}:
\begin{talign*}
        \delta\mathrm{ISS}^{(l)}  (\parameters^{(l)})  &\defeq \bar{\sigma}_{f}^{(l)}  +  \bar{\sigma}_{f}^{{(l)}^{2}} \infnorm{ R_{\tilde{h}}^{(l)} } +  \\
        &+ \frac{1}{4}  \infnorm{R_{f}^{(l)}}  \left(  \bar{\sigma}_{f}^{(l)}   \infnorm{R_{\tilde{h}}^{(l)}}   +  \bar\phi_{\tilde{h}}^{(l)}  +  1  \right), \forall l \in \layers.
\end{talign*}
Minimizing~\eqref{eq:augmented_loss} instead of~\eqref{eq:mse} aims to discourage parameter configurations that violate the sufficient condition for \dISS{}, although this formulation does not \emph{strictly} guarantee stability during the training process but rather promotes it.
Nonetheless, given that most gradient-based optimization algorithms for neural networks solve unconstrained optimization problems~\cite[Chapter 8]{goodfellowDeepLearning2016}, loss augmentation constitutes one of the most straightforward strategies that can be seamlessly implemented in any RNN training procedure.

We remark that the traditional early stopping criterion in Section~\ref{ss:training}, which is purely based on $\MSE{}(\parameters_{\kappa}; \mathcal{D}_{\mathrm{val}})$, does not account for the stability of the underlying model.
Consequently, in this work, similarly to~\cite{decarliInfinitynormbasedInputtoStateStableLong2025,bonassi_stability_2021}, we propose to update the best-found parameters $\theta^\star$ \emph{if and only if} they satisfy the \dISS{} condition in~\eqref{eq:MGU_layer_dISS_condition} for every $l \in \layers$, achieving a \emph{stability-driven early stopping} strategy.
%
%Algorithm~\ref{alg:MGU_training} details both the loss augmentation and stability-driven early stopping approaches.  

\subsection{Parameters warm-start}
\label{ss:parameters_warm_start}
As reviewed in Section~\ref{ss:training}, RNN training procedures leverage a standard initialization scheme to generate $\theta_0$.
While beneficial for general performance, this scheme initializes the MGU network parameters in a region that is \emph{likely to violate the stability conditions} derived in Section~\ref{s:stability_mgu}.
For instance, consider the ISS condition in~\eqref{eq:MGU_layer_ISS_condition}. 
The standard initialization scheme sets $\boldsymbol{b}_{f, 0}^{(l)} = \boldsymbol{1}_{n_{\hiddenunits}^{(l)}}$ for every $l \in \layers$, resulting in upper bounds of the forget gates $\bar{\sigma}_f^{(l)} \geq \sigma(1) \approx 0.73$.
With such high values for the $\bar{\sigma}_f^{(l)}$'s, the stability constraints in~\eqref{eq:MGU_layer_ISS_condition} become significantly harder to satisfy compared to an initialization characterized by a lower activation regime. 
This issue carries over to the stricter \dISS{} condition in~\eqref{eq:MGU_layer_dISS_condition}, making the standard initialization scheme likely to yield a model that violates the formal stability requirements at the start of the optimization process. 
To mitigate this, we propose a \emph{warm-start procedure} that projects the initial parameters $\theta_0$ %, obtained via standard initialization schemes, 
onto a stable \dISS{}-compliant region before MGU network training begins.
In particular, given that the \dISS{} stability condition in~\eqref{eq:MGU_layer_dISS_condition} for the $l$-th layer, $l \in \layers$, depends only on the parameters $\theta^{(l)}$ in~\eqref{eq:MGU_layer_parameters}, we consider each layer \emph{separately}.
Thus, we find the warm-start parameters $\tilde{\theta}_0^{(l)}$ by minimizing the stability violation in~\eqref{eq:augmented_loss} for the $l$-th layer:
\begin{tequation}
    \label{eq:warm_start_optimization}
    \tilde{\theta}_0^{(l)} \! \defeq \! \arg \min_{\theta^{(l)}} \left( \max \left\{ \delta\mathrm{ISS}^{(l)}(\theta^{(l)}) \! - \! (1 \! - \! \margin), 0 \right\} \right)^2 \! .
\end{tequation}
Then, MGU network training is started from $\tilde{\parameters}_0 = \{W_{y, 0}, \boldsymbol{b}_{y, 0}\}  \cup \bigcup_{l \in \layers} \tilde{\parameters}_0^{(l)}$ instead of $\theta_0$ found through the standard initialization scheme.

We suggest solving the $L$ unconstrained optimization problems in~\eqref{eq:warm_start_optimization} via \emph{local}, \emph{derivative-based}, optimization methods~\cite{NumericalOptimization2006} started at $\theta_0^{(l)}$ to attain an initial \dISS{} MGU network with parameters as close as possible to their standard initializations.
%
%Algorithm~\ref{alg:MGU_training} integrates the proposed parameters warm-start procedure within the traditional RNN training methodology.

\subsection{Projected gradient-based optimization method}
\label{ss:projected_gradient}
To \emph{strictly guarantee stability} throughout the training process, we take inspiration from the projected gradient descent scheme~\cite{boyd_convex_2023}, building on its applications in constrained RNN training~\cite{cayciConvergenceGradientDescent2025,yongGradientCentralizationNew2020}.
Let $\projection{\feasibleset}: \realset^{n_\theta} \to \feasibleset$ denote the \emph{projection operator} onto a \emph{feasible set} $\feasibleset \subset \realset^{n_\theta}$ for the parameters $\theta$ of the MGU network in~\eqref{eq:MGU_parameters}.
Then, the parameters update rule in~\eqref{eq:NN_optimizationStep} is followed by a projection onto $\feasibleset$ to ensure that $\theta_{\kappa} \in \feasibleset$ for every $\kappa \in \posintegerset$, leading to\footnote{Given that $\mathrm{GradStep}(\nabla \loss{}(\theta_{\kappa - 1}; \mathcal{D}_{\mathrm{tr}}^{(b)}), \lr)$ in~\eqref{eq:pgd_update} does not necessarily amount to $\lr \nabla \loss{}(\theta_{\kappa - 1}; \mathcal{D}_{\mathrm{tr}})$ as in the case of gradient descent, we have purposefully chosen to name our proposal \quotes{projected gradient-based optimization method} rather than \quotes{projected gradient descent} to avoid confusion.}:
\begin{tequation}
    \label{eq:pgd_update}
	\theta_{\kappa}
	=
	\projection{{\feasibleset}} \left(\theta_{\kappa - 1} \! - \! \mathrm{GradStep}\left(\nabla \loss{}(\theta_{\kappa - 1}; \mathcal{D}_{\mathrm{tr}}^{(b)}), \lr \right) \right).
\end{tequation}
Although it would be desirable to choose $\feasibleset$ as the set of parameters that satisfy the \dISS{} condition in~\eqref{eq:MGU_layer_dISS_condition} for the MGU network, this would render the projection operator $\projection{{\feasibleset}}$ \emph{computationally intractable} during RNN training due to the \emph{non-convexity} of $\feasibleset$.
We therefore focus on enforcing the simpler ISS condition in~\eqref{eq:MGU_layer_ISS_condition}.
Nonetheless, when using the projection operator in conjunction with the loss augmentation strategy in Section~\ref{ss:loss_augmentation}, guaranteeing ISS also favors the satisfaction of the \dISS{} condition in~\eqref{eq:MGU_layer_dISS_condition}.
That is because, as shown in Proposition~\ref{pr:MGU_network_dISS_imply_ISS}, the class of \dISS{} MGU networks is a subset of the class of ISS MGU networks.
Then, to ensure a computationally tractable projection, we define a stricter, \emph{convex} feasible set $\feasibleset$ as:
\begin{tequation}
    \label{eq:iss_feasible_set}
    \feasibleset \defeq \left\{ \parameters : \infnorm{ R_{\tilde{h}}^{(l)}} \leq 1 - \projectionmargin, \quad \forall l \in \layers \right\},
\end{tequation}
where $\projectionmargin \in (0, 1)$ is a small, user-defined \emph{projection margin}. %that ensures the set is closed and strictly satisfies the ISS condition.
The set in~\eqref{eq:iss_feasible_set} is stricter compared to~\eqref{eq:MGU_layer_ISS_condition} because it omits the $\bar{\sigma}_{f}^{(l)} \in (0,1)$ factor, thereby simplifying the projection to a convex operation.
In particular, the set $\feasibleset$ in~\eqref{eq:iss_feasible_set} defines an $\mathcal{L}_\infty$-norm ball constraint on the recurrent weights $R_{\tilde{h}}^{(l)}$ for every $l \in \layers$. 
Thus, the operator $\projection{\feasibleset}$ in~\eqref{eq:pgd_update} simply amounts to row-wise projections onto the $\mathcal{L}_1$-norm ball of radius $1 - \projectionmargin$, which can be computed efficiently~\cite[Chapter 8.1]{boyd_convex_2023}.
Specifically, we employ the method in~\cite{condat_fast_2016} for this purpose.
\begin{proposition}
    \label{pr:convex_projection_implies_iss}
    If $\parameters \in \feasibleset$ as in~\eqref{eq:iss_feasible_set}, then the resulting MGU network is ISS.
\end{proposition}
The proof is reported in Appendix~\ref{s:proofs}.
%
%The application of the projection operator within the MGU network training procedure is shown in Algorithm~\ref{alg:MGU_training}.

Having detailed the MGU network, its stability conditions, and the three methodologies to promote them, we now proceed to validate these approaches on numerical data.
\section{Numerical results}
\label{s:numerical_results}
In this Section, we evaluate the MGU network performance and the proposed stability-promoting training methodologies.
We first introduce the benchmark datasets and the training configuration in Section~\ref{ss:benchmark_systems}.
Then, in Section~\ref{ss:model_comparison}, we benchmark the traditional MGU network against the GRU network to validate its efficiency.
Subsequently, Section~\ref{ss:stability_comparison} evaluates the effectiveness of the proposed stability-promoting strategies.
Finally, we assess the robustness of the resulting models against noise in Section~\ref{ss:noise_robustness} and their performance on real-world data in Section~\ref{ss:silverbox}.

\subsection{Benchmark systems and datasets}
\label{ss:benchmark_systems}
We employ three distinct benchmarks to validate different properties of the MGU network.

\textbf{Dataset 1: pH Reactor.}
To compare the performance and stability properties of MGU and GRU networks, we use the pH neutralization simulator from~\cite{terzi_ph_2020}, a standard nonlinear benchmark for RNN stability~\cite{schimpernaRobustOffsetFreeConstrained2024,terzi_learning_2021}. 
The system controls the pH level ($y$) of a mixture via an alkaline flow ($u$), subject to acid and buffer disturbances.
The dataset~\cite{terzi_ph_2020} consists of $\totsequence=20$ sequences ($N^{\sequencep} = 2000, \forall v \in \{1,\ldots,V\}$) generated by a Multilevel Pseudo-Random Signal (MPRS) with sampling time $T_\mathrm{s} = 10\, \mathrm{s}$.
$\mathcal{D}$ is partitioned into $V_{\mathrm{tr}} = 15$ sequences for training, $V_{\mathrm{val}} = 4$ for validation, and $V_{\mathrm{tst}} = 1$ for testing.

\textbf{Dataset 2: Four-Tank.}
To evaluate the robustness to noise of stable MGU and GRU networks, we employ the quadruple tank benchmark system~\cite{alvaradoComparativeAnalysisDistributed2011a}, which also constitutes a standard benchmark for RNN stability evaluation~\cite{bonassi_stability_2021,schimpernaRobustConstrainedNonlinear2024}.
This system consists of four interconnected tanks receiving water from two pumps via two split valves.
It exhibits nonlinear, non-minimum phase dynamics governed by the specific coupling between the upper and lower tanks.
The task is to model the water level in the first lower tank ($y$) as a function of the two input pump flow rates ($\boldsymbol{u}$), constituting a multiple-input single-output identification problem.
Unlike the other benchmarks, where data is obtained from standard repositories, we generated this dataset by simulating the physics equations using the parameters reported in~\cite{bonassi_stability_2021}.
The system was excited with an MPRS to generate a training/validation trajectory of $2500$ steps windowed into $250$ steps subsequences ($V_{\mathrm{tr}} = 8$ and $V_{\mathrm{val}} = 2$), and an independent test trajectory ($V_{\mathrm{tst}} = 1$) of $5000$ steps, with $T_\mathrm{s} = 15\, \mathrm{s}$.
To assess robustness, additive white Gaussian noise is added to the output data at three Signal-to-Noise Ratio (SNR) levels: $100$, $10$, and $5$\footnote{The SNR is defined as the ratio between the variance of the noiseless output and the variance of the additive output noise.
We remark that prior RNN stability works, such as~\cite{terzi_learning_2021} and~\cite{bonassi_stability_2021}, involved only very low measurement noise (SNR $> 100$).}.
For each SNR level, $50$ independent training runs are performed, each characterized by a different random noise realization.
%
%Additive noise is artificially added to output data with several levels of Signal to Noise Ratio (SNR)\footnote{The SNR is defined  as the ratio between the variance of the noiseless output and the variance of the additive output noise.}.
%This controlled generation allows us to systematically inject additive white gaussian noise into the output $y$ at varying intensities, enabling a rigorous evaluation of model robustness under different measurement noise.

\textbf{Dataset 3: Silverbox.}
To validate the performance on real-world data, we utilize the Silverbox benchmark~\cite{wigrenThreeFreeData2013}.
This dataset represents a nonlinear electrical circuit simulating a mechanical vibrating system with a Duffing oscillator.
%
%It is characterized by a highly resonant and nonlinear behavior, providing a test for the modeling capacity of the MGU.
%
We use the standard data split from the literature~\cite{champneysBaselineResultsSelected2024}, dividing training/validation sequences with an 80-20\% ratio and windowing them in $250$ steps subsequences.

\textbf{Training configuration.}
As customary, we normalize inputs and outputs to the range $[-1, 1]$ to satisfy Assumption~\ref{as:input_boundedness}.
We train all networks for $e_{\mathrm{max}} = 2000$ epochs using the Adam optimizer~\cite[Chapter 8]{goodfellowDeepLearning2016} with a base learning rate of $\lr = 0.001$, decaying by a factor of $\zeta = 0.9$ every $e_{\mathrm{dc}} = 200$ epochs to aid parameters convergence.
To prevent overfitting and enhance generalization, we employ dropout with a probability $\xi = 0.05$.
Regarding the early stopping strategy, the validation mean squared error (and possibly the stability check described in Section~\ref{ss:loss_augmentation}) is computed at every epoch.
%
%For the stability-promoting methodologies, this selection is strictly subject to the satisfaction of the \dISS{} stability condition~\eqref{eq:MGU_layer_dISS_condition}\footnote{This implies that all identified models trained with stability-promoted procedures are stable.}.
%
The stability penalty weight $\isspenalty$, the stability margin $\margin$, and the projection margin $\projectionmargin$ are set to $0.01$ for all experiments.
This specific configuration was selected as it empirically demonstrated a robust balance between training convergence speed and stability enforcement across the different datasets.
The only dataset-dependent hyperparameter is the number of mini-batches $N_{\mathrm{mb}}$, which is adjusted according to the specific dataset dimensions.

All networks are implemented, trained, and validated in MATLAB 2024b on a workstation with a 64-core 3.20 GHz CPU, 2.23 GHz GPU with more than 16000 cores and 24 GB of memory, and 255 GB of RAM.
We remark that the numerical results presented here can be regarded as the \emph{first} validation of the MGU network performance on system identification tasks, as this RNN was previously applied only to more common machine learning problems such as natural language modeling, sentiment analysis, and relation extraction~\cite{zhouMinimalGatedUnit2016, liuRelationExtractionCoal2021}.

\subsection{Model comparison: MGU vs GRU networks}
\label{ss:model_comparison}
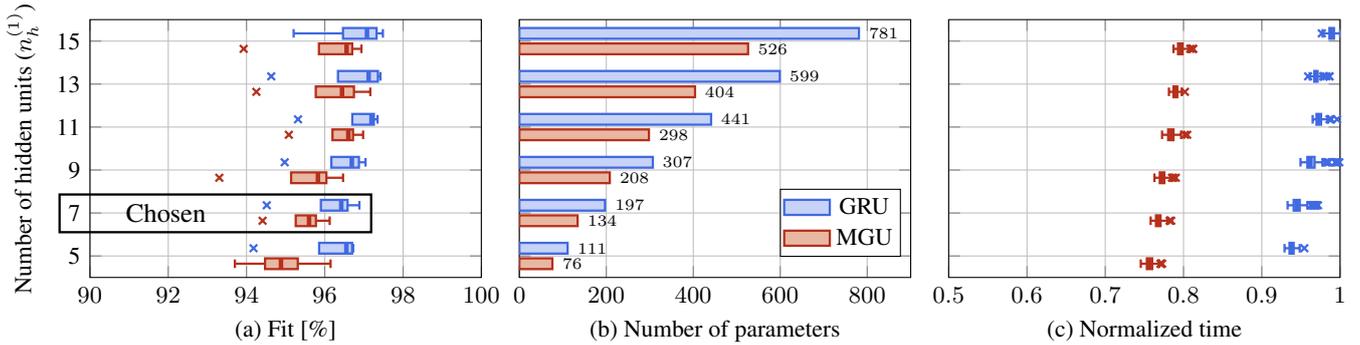
\begin{figure*}[htb]
    \centering
    % --- Parametri Grafici ---
\pgfmathsetmacro{\offsetsize}{0.18}
\newcommand{\setLow}[2]{\pgfmathsetmacro{#1}{#2-\offsetsize}}
\newcommand{\setHigh}[2]{\pgfmathsetmacro{#1}{#2+\offsetsize}}

\begin{tikzpicture}
    \begin{groupplot}[
        group style={
            group size=3 by 1,
            horizontal sep=0.5cm, % Grafici molto vicini
            vertical sep=0cm,
        },
        width=0.38\linewidth,       
        height=5cm,
        grid=major,
        % Configurazione asse Y condivisa
        ymin=0.5, ymax=6.5, % Adjusted to give space for 6
        ytick={1, 2, 3, 4, 5, 6},
        tick label style={font=\footnotesize},
        label style={font=\small},
    ]

    % ---------------------------------------------------------
    % GRAFICO 1: FIT (Boxplot)
    % ---------------------------------------------------------
    \nextgroupplot[
        xlabel={(a) $\fit{}$ [\%]},
        xmin=90, xmax=100,
        xmajorgrids=true,
        scaled x ticks=false,
        xtick distance=2,
        yticklabels={5, 7, 9, 11, 13, 15},
        ylabel={Number of hidden units ($n^{(1)}_\hiddenunits$)},
        legend image code/.code={
        \draw[#1] (0cm,-0.1cm) rectangle (0.5cm,0.1cm);
        }
    ]

    \foreach \h [count=\i] in {5, 7, 9, 11, 13, 15} {
    \pgfmathsetmacro{\colMgu}{int((\i-1)*2)}
    \pgfmathsetmacro{\colGru}{int((\i-1)*2+1)}
    \setLow{\posMgu}{\i}
    \setHigh{\posGru}{\i}

    % IMPORTANTE: forget plot qui impedisce la creazione di legende automatiche sbagliate
    \addplot+ [mgu box style, boxplot/draw position=\posMgu, forget plot] 
        table [y index=\colMgu] {data/mgu_gru_hu.dat};

    \addplot+ [gru box style, boxplot/draw position=\posGru, forget plot] 
        table [y index=\colGru] {data/mgu_gru_hu.dat};
    }

    % ---------------------------------------------------------
    % GRAFICO 2: PARAMETERS
    % ---------------------------------------------------------
    \nextgroupplot[
        xlabel={(b) Number of parameters},
        yticklabels={},         
        ylabel={},
        xmin=0, xmax=900,
        xmajorgrids=true,
        nodes near coords style={
            font=\tiny, 
            color=black, 
            anchor=west, 
            /pgf/number format/.cd, fixed, precision=0, 1000 sep={} % Formattazione intera pulita senza virgole
        },
        legend style={
            at={(0.975, 0.08)},
            anchor=south east,
            legend columns=1,
            draw=black,
            font=\small,
            /tikz/every even column/.append style={column sep=0.3cm},
            row sep=0cm,
            inner sep=1pt,
        },
    ]

    % --- LEGENDA MANUALE ---
    % Usiamo gli stili 'legend' che non hanno la logica boxplot    
    \addlegendimage{area legend, gru legend}
    \addlegendentry{GRU}
    
    \addlegendimage{area legend, mgu legend}
    \addlegendentry{MGU}
    
    % MGU Bars
    \addplot+ [mgu bar style, point meta=explicit] 
    table [x=mgu_params, y expr=\coordindex+1-\offsetsize, meta=mgu_params] {model_params.dat};

    % GRU Bars
    \addplot+ [gru bar style, point meta=explicit] 
    table [x=gru_params, y expr=\coordindex+1+\offsetsize, meta=gru_params] {model_params.dat};

    % ---------------------------------------------------------
    % GRAFICO 3: TIME (Boxplot)
    % ---------------------------------------------------------
    \nextgroupplot[
        xmin=0.5, xmax=1,
        xlabel={(c) Normalized time},
        xtick distance=0.1,
        yticklabels={},         
        ylabel={},
        legend image code/.code={
             \draw[#1] (0cm,-0.1cm) rectangle (0.5cm,0.1cm);
        }
    ]

    \foreach \h [count=\i] in {5, 7, 9, 11, 13, 15} {
    \pgfmathsetmacro{\colMgu}{int((\i-1)*2)}
    \pgfmathsetmacro{\colGru}{int((\i-1)*2+1)}
    \setLow{\posMgu}{\i}
    \setHigh{\posGru}{\i}

    % IMPORTANTE: forget plot qui impedisce la creazione di legende automatiche sbagliate
    \addplot+ [mgu box style, boxplot/draw position=\posMgu, forget plot] 
        table [y index=\colMgu] {data/mgu_gru_times.dat};

    \addplot+ [gru box style, boxplot/draw position=\posGru, forget plot] 
        table [y index=\colGru] {data/mgu_gru_times.dat};
    }
    \end{groupplot}

    % Draw the red framing box with a text on top of it
    \draw [
        line width=\linewidthlight,
        draw=black,
    ] (-0.4, 0.6) rectangle (3.7, 1.1);

    \node at (1, 0.85) [
            anchor=center,
            font=\color{black}
        ] {Chosen};
\end{tikzpicture}
    \vspace{-0.7cm}
    \caption{Comparison analysis between MGU and GRU networks for varying numbers of hidden units ($n^{(1)}_\hiddenunits{}$) on Dataset 1 (pH Reactor).
    (a) Distribution of the $\fit{}$ in~\eqref{eq:fit} over the different sequences that compose $\mathcal{D}_{\mathrm{val}}$, demonstrating comparable performance between MGU and GRU networks.
    (b) Parameters count, showing that MGU networks require approximately two-thirds of the GRU network parameters (see Remark~\ref{rem:mgu_smaller}).
    (c) Distribution of the normalized inference times, illustrating the computational efficiency of the MGU network architecture.
    From this analysis, we select $n^{(1)}_\hiddenunits{} = 7$ for subsequent stability experiments on Dataset 1.
    }
    \label{fig:rnn_comparison}
\end{figure*}
We first establish the baseline performance of the MGU network against the widely used GRU network using the traditional training methodology reviewed in Section~\ref{ss:training} on Dataset 1 (Section~\ref{ss:benchmark_systems}).
The optimal model complexity is determined by training a single-layer ($L = 1$) network, varying $n^{(1)}_\hiddenunits{}$ from 5 to 15 (in increments of 2), minimizing only the $\MSE{}(\parameters;\mathcal{D}_{\mathrm{tr}})$ in~\eqref{eq:mse}.
We evaluate the performance on the different sequences that compose the validation dataset $\mathcal{D}_{\mathrm{val}}$ using the $\fit{}$ metric:
\begin{tequation}
    \label{eq:fit}
    \fit{}^{\sequencep}\left(\parameters^\star\right) \defeq{} 100 \left( 1 - \frac{\sqrt{\sum_{k = 0}^{N^{\sequencep} - 1}\left( y_{k}^{\sequencep} - \hat{y}_{k}^{\sequencep} \left(\parameters^\star\right) \right)^2}}{\sqrt{\sum_{k = 0}^{N^{\sequencep} - 1}\left( y_{k}^{\sequencep} - \bar{y}^{\sequencep} \right)^2}} \right),
\end{tequation}
where $\theta^\star$ denotes the identified parameters (see Algorithm~\ref{alg:MGU_training}).
The results in \figurename{}~\ref{fig:rnn_comparison} highlight three findings.
First, Panel~(a) shows that \emph{MGU networks achieve predictive accuracies comparable to GRU networks} across all network sizes.
Second, Panel~(b) confirms that \emph{an MGU network requires approximately two-thirds of the parameters of an equivalent GRU network}, consistent with Remark~\ref{rem:mgu_smaller}.
Third, this structural efficiency translates to a \emph{reduction in computational cost}: Panel~(c) shows lower normalized inference times per step\footnote{Times are normalized by the slowest result observed, which is in the order of $0.1\,\mathrm{ms}$. %using MATLAB 2024b.
}, and training time was reduced from 1.7 hours (GRU) to 1.2 hours (MGU)\footnote{Based on the MATLAB 2024b custom layer implementation.}.
Based on these results, we select $n_\hiddenunits^{(1)} = 7$ for all subsequent experiments on Dataset 1, matching the complexity used in~\cite{terzi_learning_2021} where stable LSTM networks were employed for the same dataset.

\subsection{Stability comparison}
\label{ss:stability_comparison}
\begin{figure*}
    \centering
    \tikzexternaldisable
    % --- Parametri Grafici ---
\pgfmathsetmacro{\offsetsize}{0.18}
\newcommand{\setLow}[2]{\pgfmathsetmacro{#1}{#2-\offsetsize}}
\newcommand{\setHigh}[2]{\pgfmathsetmacro{#1}{#2+\offsetsize}}
\def\verticalLinePercent{96.5} % for dISS LSTM comparison

\begin{tikzpicture}
    \begin{groupplot}[
        group style={
            group size=2 by 1,
            horizontal sep=0.9cm,
            vertical sep=0cm,
        },
        width=0.435\linewidth, 
        height=4cm,
        grid=major,
        % Shared Y axis settings (identical in both files)
        ymin=0, ymax=6,
        ytick={1, 2, 3, 4, 5},
        tick label style={font=\footnotesize},
        label style={font=\small},
    ]
    
    % ---------------------------------------------------------
    % GRAFICO 1: FIT (Boxplot)
    % ---------------------------------------------------------
    \nextgroupplot[
        xlabel={(a) $\fit{}$ [\%]},
        xmin=0, xmax=100,
        xmajorgrids=true,
        scaled x ticks=false,
        xtick distance=10,
        ylabel={Training method},
        yticklabels={MSE, LA, WS, PGM, PGM+WS},
        legend style={
        at={(rel axis cs:0.01,0.98)},%,0.635)},
        anchor=north west,
        draw=black
      }
    ]   
        % none
        \setLow{\posMgu}{1}
        \setHigh{\posGru}{1}
        \addplot+ [mgu box style, boxplot/draw position=\posMgu, forget plot] table [y=fit_none] {data/performance_data.dat};
        \addplot+ [gru box style, boxplot/draw position=\posGru, forget plot] table [y=fit_none_gru] {data/performance_data.dat};
        
        % aug
        \setLow{\posMgu}{2},
        \setHigh{\posGru}{2},
        \addplot+ [mgu box style, boxplot/draw position=\posMgu, forget plot] table [y=fit_aug] {data/performance_data.dat};
        \addplot+ [gru box style, boxplot/draw position=\posGru, forget plot] table [y=fit_aug_gru] {data/performance_data.dat};
        
        % aug + init
        \addplot+ [mgu box style, boxplot/draw position=3, forget plot] table [y=fit_aug_init] {data/performance_data.dat};
        
        % aug + proj
        \addplot+ [mgu box style, boxplot/draw position=4, forget plot] table [y=fit_aug_proj] {data/performance_data.dat};
        
        % all
        \addplot+ [mgu box style, boxplot/draw position=5, forget plot] table [y=fit_aug_proj_init] {data/performance_data.dat};   

        \addplot [dashed, line width={\linewidthlight}, black] coordinates {(\verticalLinePercent, 0) (\verticalLinePercent, 7)};
        \addlegendentry{LSTM$_\text{LA}$~\cite{terzi_learning_2021}}

        % Legend by color
        \addlegendimage{area legend, gru legend}
        \addlegendentry{GRU}
    
        \addlegendimage{area legend, mgu legend}
        \addlegendentry{MGU}
    
    % ---------------------------------------------------------
    % GRAFICO 3: dIUSS in-range
    % ---------------------------------------------------------
    \nextgroupplot[
        xmin=0, xmax=100,
        xlabel={(b) \dISS{} in-range rate $[\%]$},
        xtick distance=10,
        yticklabels={},         
        ylabel={},
        %trim axis left, trim axis right,
    ]   
        % none
        \setLow{\posMgu}{1},
        \setHigh{\posGru}{1},
        \addplot+ [mgu box style, boxplot/draw position=\posMgu, forget plot] table [y=dISS_none] {data/performance_data.dat};
        \addplot+ [gru box style, boxplot/draw position=\posGru, forget plot] table [y=dISS_none_gru] {data/performance_data.dat};
        
        % aug
        \setLow{\posMgu}{2},
        \setHigh{\posGru}{2},
        \addplot+ [mgu box style, boxplot/draw position=\posMgu, forget plot] table [y=dISS_aug] {data/performance_data.dat};
        \addplot+ [gru box style, boxplot/draw position=\posGru, forget plot] table [y=dISS_aug_gru] {data/performance_data.dat};
        
        % aug + init
        \addplot+ [mgu box style,    boxplot/draw position=3, forget plot] table [y=dISS_aug_init] {data/performance_data.dat};
        
        % aug + proj
        \addplot+ [mgu box style,  boxplot/draw position=4, forget plot] table [y=dISS_aug_proj] {data/performance_data.dat};
        
        % all
        \addplot+ [mgu box style,boxplot/draw position=5, forget plot] table [y=dISS_aug_proj_init] {data/performance_data.dat};   
  \end{groupplot}

    % Box
    \draw [
        line width=\linewidthlight,
        draw=black,
    ] (5, 1.0) rectangle (12.9, 1.4);

    % Text 1
    \node at (7.6, 1.2) [
            anchor=west,
            font=\color{black}
        ] {Preferred method: WS};

    % Text 2
    \node at (7.6, 0.4) [
            anchor=west,
            font=\color{black}
        ] {Never in range};
    \draw [->, >=stealth, thick] (7.6, 0.4) -- (7.2, 0.4);
        
\end{tikzpicture}
    \tikzexternalenable
    \vspace{-0.4cm}
    \caption{Comparison of stability-promoting methods over 30 training runs on Dataset 1 (pH Reactor), highlighting the performance-stability trade-off.
    (a) Distribution of the $\fit{}$ in~\eqref{eq:fit} on the test set $\mathcal{D}_{\mathrm{tst}}$.
    (b) Distribution of the \dISS{} in-range rate, showing the percentage of epochs satisfying the \dISS{} stability condition in~\eqref{eq:MGU_layer_dISS_condition} for MGU models and the corresponding condition from~\cite{bonassi_stability_2021} for GRU models.
    Note the balance achieved by MGU$_{\text{WS}}$ and the low-accuracy of PGM and PGM+WS methodologies.
    }
    \label{fig:methods_comparison}
\end{figure*}
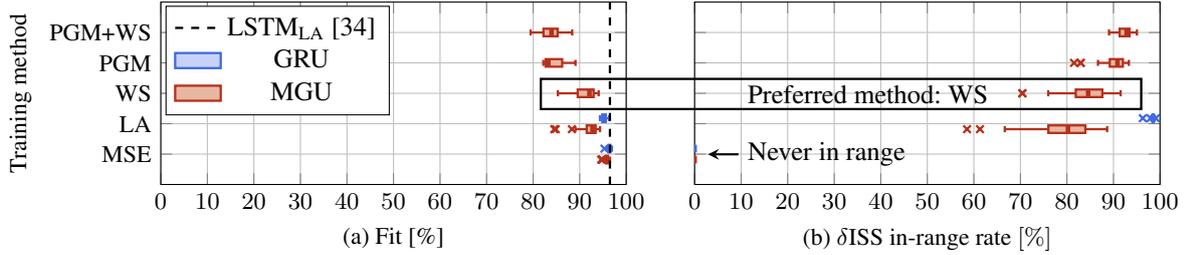
We now evaluate the three stability-promoting strategies introduced in Section~\ref{s:stability_promoted_training}: Loss Augmentation ($\text{LA}$), parameters Warm-Start ($\text{WS}$), and Projected Gradient-based optimization Method ($\text{PGM}$).
To ensure clarity, we adopt the nomenclature $\text{RNN}_{\text{METHOD}}$ for our experimental configurations, where RNN identifies the network architecture (either MGU or GRU), and the subscript denotes the training methodology:
\begin{itemize}
    \item $\text{RNN}_{\text{MSE}}$: traditional training with only the MSE as the loss function (Section~\ref{ss:training}).
    \item $\text{RNN}_{\text{LA}}$: training with loss augmentation (Section~\ref{ss:loss_augmentation}).
    \item $\text{RNN}_{\text{WS}}$: training with LA and warm-start (Section~\ref{ss:parameters_warm_start}).
    \item $\text{RNN}_{\text{PGM}}$: training with LA and PGM (Section~\ref{ss:projected_gradient}).
    \item $\text{RNN}_{\text{PGM+WS}}$: training with LA, WS, and PGM.
\end{itemize}
We remark that $\text{RNN}_{\text{MSE}}$ leverages the traditional early stopping criterion purely based on $\MSE{}(\parameters; \mathcal{D}_{\mathrm{val}})$, while the other four methodologies employ the stability-driven early stopping strategy introduced in Section~\ref{ss:loss_augmentation}.
In this comparison, we consider only single-layer architectures for all models.
While MGU networks are tested with all the proposed methods, GRU networks are trained only with MSE and LA, as these represent the standard benchmarks available in the literature for GRU networks~\cite{bonassi_stability_2021}.
%
%We train GRU networks only with MSE and LA, while testing MGU networks with all the proposed methods.
%
We set $n^{(1)}_\hiddenunits{} = 7$ according to Section~\ref{ss:model_comparison}, and evaluate the performance over $30$ independent training runs, each involving different initial parameters $\parameters_{0}$ due to the \emph{randomness} introduced by the Glorot and orthogonal initializations (Section~\ref{ss:training}), shared across the various configurations.
When using parameters warm-start, we solved the optimization problems in~\eqref{eq:warm_start_optimization} with a quasi-Newton method~\cite[Chapter 6]{NumericalOptimization2006}. % with high precision tolerances.

\figurename{}~\ref{fig:methods_comparison} illustrates the results, highlighting the trade-off between predictive accuracy and stability promotion.
Panel~(a) reports the $\fit{}$ distribution in~\eqref{eq:fit} on $\mathcal{D}_{\mathrm{tst}}$, while Panel~(b) depicts the \dISS{} in-range rate, defined as the percentage of training epochs where the parameters satisfy the \dISS{} condition in~\eqref{eq:MGU_layer_dISS_condition}.
Although the stability-driven early stopping criterion (Section~\ref{ss:loss_augmentation}) ensures that the final parameters $\parameters^\star$ belong to the \dISS{}-compliant region, the \dISS{} in-range rate serves as a metric for the optimization landscape's tractability and the method's efficiency in maintaining the search within the \dISS{}-compliant region.

Regarding predictive performance, we observe that all methods achieve a satisfactory median $\fit{}$.
The unconstrained $\text{MGU}_{\text{MSE}}$ and GRU$_{\text{MSE}}$ models exhibit the highest median $\fit{}$ with the lowest variance.
As expected, \emph{introducing stability constraints restricts the feasible parameter space, resulting in a slightly lower median $\fit{}$ and a modest increase in variance}, a common trade-off when imposing penalties or hard constraints within an optimization problem~\cite[Chapter 4.7.5]{boyd_convex_2023}.
For context, the single-layer LSTM$_{\text{LA}}$ model ($n^{(1)}_\hiddenunits{} = 7$) from~\cite{terzi_learning_2021} achieved a $\fit{}$ of 96.5\%.
MGU$_{\text{LA}}$ and GRU$_{\text{LA}}$ models perform slightly below this benchmark, a discrepancy likely attributable to the LSTM network's dual-state architecture and higher parameter count, or favorable initial parameter set~\cite{murphyProbabilisticMachineLearning2025}.

\begin{figure*}
    \centering
    % Points to plot each for compactness
\def\npoints{5}
\pgfplotsset{
    % --- Styles for lines/fills (unchanged) ---
    real/.style={
        draw=black,
        very thick,
        densely dashed,
        smooth
    },
    mgu_med/.style={
        draw=mgu_color,
        very thick,
        smooth
    },
    gru_med/.style={
        draw=gru_color,
        very thick,
        smooth
    },
    mgu_fill/.style={
        fill=mgu_color,
        fill opacity=0.2
    },
    gru_fill/.style={
        fill=gru_color,
        fill opacity=0.2
    },
}
% --- Define the plotting command ---
% #1 = data name prefix (e.g., "none")
% #2 = style prefix (e.g., "mgu", "gru")
% #3 = 'true' to plot legend entries, otherwise 'false' (or anything else)
\newcommand{\generaterealizationplot}[3]{
    % Plot the upper bound
    \addplot+[name path=max, #2_color, forget plot] table [x=k, y=#1_max] {data/realization_data.dat};
    % Plot the lower bound
    \addplot+[name path=min, #2_color, forget plot] table [x=k, y=#1_min] {data/realization_data.dat};
    % Fill between the paths (Uses forget plot to prevent automatic legend entry)
    \addplot[#2_fill, forget plot] fill between [of=max and min];
    % Add Min-Max Legend Entry
    \ifstrequal{#3}{true}{
        % Overwrite the legend image for Min-Max with the custom filled box
        \addlegendimage{
            legend image code/.code={
                \draw[draw=#2_color, fill=#2_color!30, line width=0.4pt] (0cm,-0.08cm) rectangle (0.4cm,0.08cm);
            }
        }
    \addlegendentry{Min-Max}
    }{}
    % Plot the median line on top 
    \addplot+[#2_med] table [x=k, y=#1_med, name path=median] {data/realization_data.dat};
    % Add Median Legend Entry
    \ifstrequal{#3}{true}{
    \addlegendentry{Median}
    }{}
    % Plot the real (dashed) line on top (always)
    \addplot+[real] table [x=k, y=Y_real, name path=true] {data/realization_data.dat};
    % Add Real Legend Entry
    \ifstrequal{#3}{true}{
    \addlegendentry{Real}
    }{}
}
% --- Add some spacing between the table rows ---
\renewcommand{\arraystretch}{1.5}
\begin{tikzpicture}
\begin{groupplot}[
    % --- Configure the grid ---
    group style={
        group size=2 by 3, % 2 columns, 3 rows
        vertical sep=0.1cm,  % Space between rows
        horizontal sep=1.5cm   % Space between columns
    },
    % --- Common styles for ALL plots ---
    width=0.9\columnwidth, % Use slightly less than 0.5 for safety
    height=3.5cm,
    grid=major,
    no markers,
    xmin=0, xmax=20000,
    ymin=5.5, ymax=8.5,
    % ticks 
    ytick={6.0, 6.5, ..., 8.0},
    xtick distance=2000,
    %x tick label style={font=\scriptsize, /pgf/number format/fixed,},  
    %y tick label style={font=\scriptsize, /pgf/number format/fixed,},  
    scaled x ticks=false,
    each nth point={\npoints},
    legend style={
        legend columns=3,
        %/pgf/number format/legend identifier width=6pt, % Riduci la larghezza orizzontale (es. da 12pt a 6pt)
        %/tikz/every even column/.append style={column sep=2pt}, % Riduci lo spazio tra le colonne della legenda
        %font=\scriptsize, % Manteniamo il font piccolo per coerenza
        nodes={scale=0.8}, 
        inner sep=1pt,
    }
    ]
    % --- Row 1 ---
    \nextgroupplot[ylabel={$\text{GRU}_{\text{MSE}}$}, xticklabels=\empty, font=\scriptsize]
    \generaterealizationplot{gru_none}{gru}{false}
    \nextgroupplot[ylabel={$\text{GRU}_{\text{LA}}$}, xticklabels=\empty, font=\scriptsize]
    \generaterealizationplot{gru_aug}{gru}{true}
    % --- Row 2 ---
    \nextgroupplot[ylabel={$\text{MGU}_{\text{MSE}}$}, xticklabels=\empty, font=\scriptsize]
    \generaterealizationplot{none}{mgu}{false}
    \nextgroupplot[ylabel={$\text{MGU}_{\text{LA}}$}, xticklabels=\empty, font=\scriptsize]
    \generaterealizationplot{aug}{mgu}{true}
    % --- Row 3 ---
    \nextgroupplot[ylabel={$\text{MGU}_{\text{WS}}$}, 
    scaled x ticks=true,
    xlabel={Time [s]},
    font=\scriptsize
    ]
    \generaterealizationplot{aug_init}{mgu}{false}
    \nextgroupplot[
     ylabel={$\text{MGU}_{\text{PGM+WS}}$},
    scaled x ticks=true,
    xlabel={Time [s]},
    font=\scriptsize
    ]
    \generaterealizationplot{aug_proj_init}{mgu}{false}
    %\nextgroupplot[
    % ylabel={$\text{MGU}_{\text{PGD}}$},
    %scaled x ticks=true,
    %xlabel={Time [s]},
    %]
    %\generaterealizationplot{aug_proj}{mgu}{false}
\end{groupplot}
\end{tikzpicture}
    \vspace{-0.3cm}
    \caption{Time-domain output on $\mathcal{D}_{\mathrm{tst}}$ of Dataset 1 (pH Reactor) for the median-performing model (solid lines) from each configuration.
    The filled region indicates the Min-Max range across the 30 training runs.
    The ground truth (real system output) is shown as a dashed black line.
    }
     \label{fig:realization_comparison}
\end{figure*}
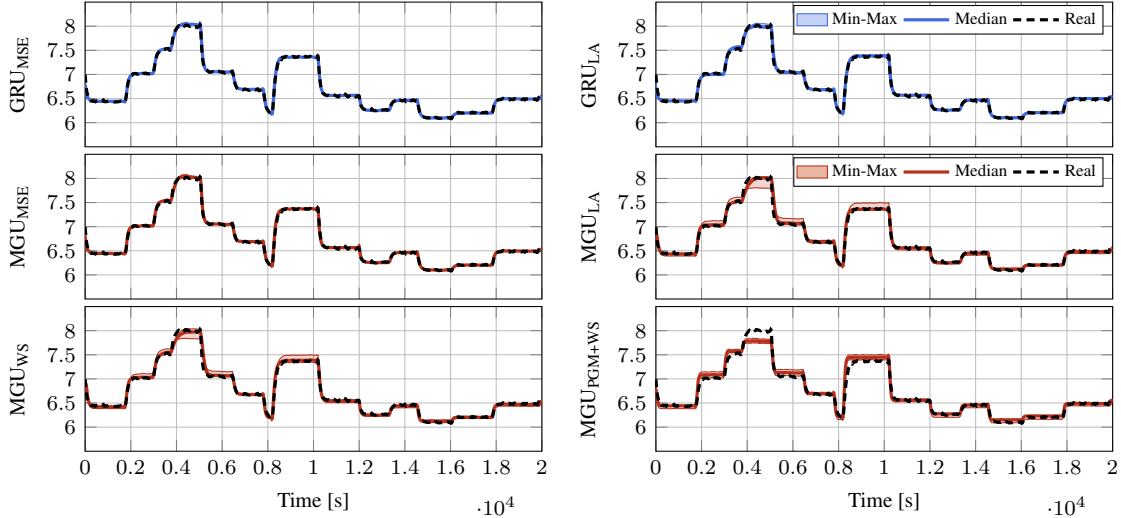

In terms of stability, the $\text{MGU}_{\text{MSE}}$ and GRU$_{\text{MSE}}$ models \emph{never} enter the \dISS{}-compliant region, confirming that \emph{unconstrained training does not inherently yield \dISS{} models}.
Conversely, the GRU$_{\text{LA}}$ model exhibits a high median \dISS{} in-range rate, which we attribute to the $\text{GRU}$ network's large parameter count, making its stability constraint less restrictive~\cite[Chapter 8]{goodfellowDeepLearning2016}.
For the MGU network, we observe a distinct progression in stability adherence as additional strategies are employed.
While MGU$_{\text{LA}}$ achieves a median \dISS{} in-range rate of 80\%, the inclusion of the parameters warm-start strategy (MGU$_{\text{WS}}$) increases this to 85\%.
Notably, the computational overhead for this warm-start was negligible compared to the total training duration, requiring a few seconds to solve the optimization problems in~\eqref{eq:warm_start_optimization}.
The PGM-based methods show the strongest adherence to stability, with median \dISS{} in-range rates exceeding 90\%.
However, this comes at the cost of a significant drop in predictive performance, suggesting that the projection operation may excessively limit the parameter search space required for accurate modeling.
Considering the trade-off between model $\fit{}$ and the consistency of stability adherence during training, \emph{MGU$_{\text{WS}}$ offers the most effective balance}, and it will be used for the subsequent experiments.
%
%It provides substantial stability \dISS{} in-range improvement over MGU$_{\text{LA}}$ while maintaining a high median $\fit{}$, confirming parameters warm-start as our preferred method, used for the subsequent experiments.
%
Finally, in \figurename{}~\ref{fig:realization_comparison}, we display the time-domain output for the median-performing model from each configuration.
Consistent with the performance distributions shown in \figurename{}~\ref{fig:methods_comparison}, the GRU configurations and the MGU$_{\text{MSE}}$ network appear relatively robust to the choice of the initial parameters $\parameters_{0}$.
In contrast, \emph{the stability-promoted MGU models exhibit a higher sensitivity to the initialization}.
Finally, while the median predictions of almost all configurations closely track the ground-truth data (black dashed line), the MGU$_{\text{PGM+WS}}$ displays a less accurate tracking of steady-state pH levels and a visibly lower $\fit{}$ during large disturbances.
%
%The medians for the GRU and MGU standard MSE and LA configurations are visually almost identical and closely track the ground-truth data (black dashed line), though the MGU$_{\text{LA}}$ shows a higher variability range.
%
%The PGD-based methods (MGU$_{\text{PGD}}$ and MGU$_{\text{PGD+WS}}$), on the other hand, show some minor inaccuracies in modeling, with a visibly worse $\fit{}$ during large disturbances.
%

\subsection{Robustness to noise on output measurements}
\label{ss:noise_robustness}
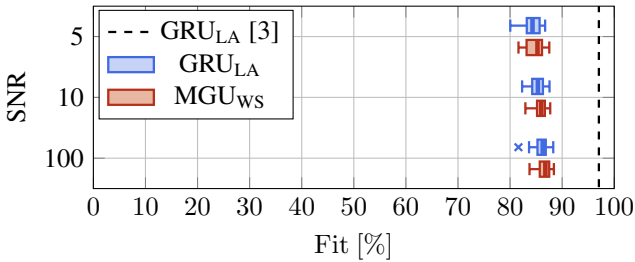
\begin{figure}
    %\centering
    \tikzexternaldisable
    % --- Parametri Grafici ---
\pgfmathsetmacro{\offsetsize}{0.18}
\newcommand{\setLow}[2]{\pgfmathsetmacro{#1}{#2-\offsetsize}}
\newcommand{\setHigh}[2]{\pgfmathsetmacro{#1}{#2+\offsetsize}}
\def\verticalLinePercent{97.05} % for dISS GRU Bonassi comparison

\begin{tikzpicture}
    \begin{axis}[
        width=\columnwidth,
        height=4cm,
        ymin=0.5, ymax=3.5,
        ytick={1, 2, 3},
        yticklabels={100, 10, 5},
        grid=major,
        xmin=0, xmax=100,
        ylabel={SNR},
        xlabel={$\mathrm{Fit}~[\%]$},
        xtick distance=10,
        legend style={
            at={(rel axis cs:0.01,0.98)},%,0.635)},
            anchor=north west,
            draw=black
        }
    ]
    % --- PLOT DATA (Fit) ---
    
    % 100
    \setLow{\posMgu}{1}
    \setHigh{\posGru}{1}
    \addplot+ [mgu box style, boxplot/draw position=\posMgu, forget plot] table [y=Fit_mgu_WS_SNR100] {data/performance_snr.dat};
    \addplot+ [gru box style, boxplot/draw position=\posGru, forget plot] table [y=Fit_gru_LA_SNR100] {data/performance_snr.dat};
    
    % 10
    \setLow{\posMgu}{2}
    \setHigh{\posGru}{2}
    \addplot+ [mgu box style, boxplot/draw position=\posMgu, forget plot] table [y=Fit_mgu_WS_SNR10] {data/performance_snr.dat};
    \addplot+ [gru box style, boxplot/draw position=\posGru, forget plot] table [y=Fit_gru_LA_SNR10] {data/performance_snr.dat};
    
    % 5
    \setLow{\posMgu}{3}
    \setHigh{\posGru}{3}
    \addplot+ [mgu box style, boxplot/draw position=\posMgu, forget plot] table [y=Fit_mgu_WS_SNR5] {data/performance_snr.dat};
    \addplot+ [gru box style, boxplot/draw position=\posGru, forget plot] table [y=Fit_gru_LA_SNR5] {data/performance_snr.dat};

    % --- Comparison Line (GRU bonassi) ---
    \addplot [dashed, line width={\linewidthlight}, black] coordinates {(\verticalLinePercent, 0) (\verticalLinePercent, 7)};
    
    \addlegendentry{GRU$_\text{LA}$~\cite{bonassi_stability_2021}}

    % Legend by color
    \addlegendimage{area legend, gru legend}
    \addlegendentry{GRU$_\text{LA}$}
    
    \addlegendimage{area legend, mgu legend}
    \addlegendentry{MGU$_\text{WS}$}
    
    \end{axis}
\end{tikzpicture}
    \tikzexternalenable
    \vspace{-0.7cm}
    \caption{Performance comparison between the MGU$_\text{WS}$ and GRU$_\text{LA}$ models on Dataset 2 (Four-Tank) under varying SNRs.
    The distributions are computed over 50 independent runs with different random noise realizations.
    Both RNNs maintain high predictive accuracy even under severe noise (SNR=5).
    The dashed line indicates the performance of a GRU$_\text{LA}$ model trained with low measurement noise (taken from~\cite{bonassi_stability_2021}).
    }
    \label{fig:snr_comparison}
\end{figure}
In this Section, we evaluate the robustness to measurement noise of stable RNNs using Dataset 2 (Section~\ref{ss:benchmark_systems}).
Particularly, the MGU$_{\text{WS}}$ model is compared against the GRU$_{\text{LA}}$ one.
We employ a 3-layer architecture with $n_{\hiddenunits}^{(l)}=7$ hidden units in each layer $l \in \layers$, matching the $\text{GRU}_{\text{LA}}$ model used in the literature benchmark~\cite{bonassi_stability_2021}.
A washout period of $25$ steps is applied to the training and validation subsequences. 
This means that the prediction error during these initial steps is excluded from the $\MSE{}$ calculation in~\eqref{eq:mse} to mitigate the influence on the loss function of the initial transient due to a potentially wrong initial state of the network~\cite{goodfellowDeepLearning2016}.
The results are summarized in \figurename{}~\ref{fig:snr_comparison}. 
The MGU$_{\text{WS}}$ model maintains a predictive accuracy that is competitive with, and occasionally superior to, that of the GRU$_{\text{LA}}$ model across all noise levels, \emph{despite having fewer parameters}.
%
%Specifically, at the lowest SNR of $5$ the median performance of the stable RNNs remains robust, demonstrating that these models not degrade their ability to generalize in noisy conditions.
%
For comparison with available results from the literature, the GRU$_{\text{LA}}$ model reported in~\cite{bonassi_stability_2021}, trained on the same system but this time with a low measurement noise (SNR $> 100$), achieves higher performance (\fit{} of 97.05\%).

\subsection{Silverbox performance}
\label{ss:silverbox}
Finally, we compare the MGU$_{\text{WS}}$ and GRU$_{\text{LA}}$ models on real-world data using Dataset 3 (Section~\ref{ss:benchmark_systems}).
To identify the optimal model structure, we conducted a grid search inspired by the protocol in~\cite{champneysBaselineResultsSelected2024}.
Specifically, we explored architectures with depth $L \in \{1, 2, 3\}$ and hidden units $n_{\hiddenunits}^{(l)} \in \{8, 16, 32, 64\}$ per layer.
Following the methodology in~\cite{champneysBaselineResultsSelected2024}, to ensure a fair comparison with literature baselines, we tested different input lag windows $n_{\text{lags}} \in \{5, 10, 15, 20\}$. 
Specifically, the network is fed at each time step $k$ with an augmented input vector $\boldsymbol{u}_{k,\text{lags}} = [u_k, \dots, u_{k-n_{\text{lags}}+1}]^\top \in \realset^{n_{\text{lags}} n_u}$, where $n_u=1$ is the dimensionality of the original input.
Consistent with Section~\ref{ss:noise_robustness}, we use a washout period of $25$ steps to mitigate the effect of initial transients on the loss function.
Based on the lowest $\MSE{}$ recorded on $\mathcal{D}_{\mathrm{val}}$, the optimal configuration for both networks utilizes an input lag of $n_{\text{lags}}=10$ and a depth of $L=3$ layers.
Further, the MGU$_{\text{WS}}$ model uses $n_{\hiddenunits}^{(l)}=64$ units per layer, while the GRU$_{\text{LA}}$ network $n_{\hiddenunits}^{(l)}=8$, $\forall l \in \layers$.

We report the performance on the three standard test configurations for the Silverbox benchmark: the \textit{Multisine} sequence and the \textit{Arrow} trajectory, evaluated both on the \textit{Full} sequence, and the \textit{No Extrapolation} subset~\cite{wigrenThreeFreeData2013,champneysBaselineResultsSelected2024}.
The evaluation metric is the square root of the MSE (RMSE).
The results are summarized in Table~\ref{tab:silverbox_results}.
The MGU$_{\text{WS}}$ model shows \emph{superior performance}, outperforming the baseline GRU$_{\text{LA}}$ network by a margin across all test configurations.
We emphasize that the GRU$_{\text{LA}}$ model results, despite being suboptimal, represent the best performance obtained across the entire hyperparameter grid, as also other configurations consistently failed to yield accurate models.
While the GRU$_{\text{LA}}$ network struggled to capture the dynamics of the system, 
%likely getting trapped in a poor local minimum due to the non-convexity of the loss landscape under stability constraints, 
the MGU$_{\text{WS}}$ strategy successfully converged to a satisfactory model.

We also benchmark our results against \emph{unconstrained} baselines (i.e., without stability concerns) reported in the recent survey~\cite{champneysBaselineResultsSelected2024}.
On the Arrow (Full) test set, the standard unconstrained GRU network from~\cite{champneysBaselineResultsSelected2024} achieved an RMSE of $2.80\,\mathrm{mV}$, while the simple RNN achieved an RMSE of $4.88\,\mathrm{mV}$.
Our stable MGU$_{\text{WS}}$ model achieves an RMSE of $6.88\,\mathrm{mV}$.
While this represents an increase in error compared to the unconstrained state-of-the-art, it is quantitatively comparable to the standard RNN baseline and other methods in~\cite{champneysBaselineResultsSelected2024}.
This result highlights the \emph{accuracy-stability trade-off}: the proposed method sacrifices a marginal amount of predictive accuracy to ensure rigorous \dISS{} stability guarantees.
\begin{table}[t]
    \centering
    \caption{Performance comparison on Dataset 3 (Silverbox) with an MGU$_{\text{WS}}$ network of $n_{\hiddenunits}^{(l)}=64$ units per layer and a GRU$_{\text{LA}}$ model of $n_{\hiddenunits}^{(l)}=8$ for all $l \in \layers$, both considering $n_{\text{lags}}=10$ lags and with 3 layers. 
    The columns correspond to the standard test splits defined in~\cite{champneysBaselineResultsSelected2024}. 
    Results are reported in RMSE ($\mathrm{mV}$).
    Best results are reported with a bold font.
    }
    \begin{tabular}{lccc}
        \textbf{Model} & \textbf{Multisine} & \textbf{Arrow (Full)} & \textbf{Arrow (No Extrap.)} \\ \hline
        \colorrowoftable MGU$_{\text{WS}}$ & \textbf{5.27} & \textbf{6.88} & \textbf{4.37} \\
         GRU$_{\text{LA}}$ & 24.67 & 24.38 & 19.37\\ \hline
    \end{tabular}
    \label{tab:silverbox_results}
\end{table}
\section{Conclusion}
\label{s:conclusion}
This work investigated the use of MGU networks for the identification of stable nonlinear models, with a particular focus on deriving sufficient ISS and \dISS{} parametric conditions for the MGU model. 
To promote compliance with the \dISS{} condition, we also proposed and compared three distinct training methodologies: loss augmentation, parameters warm-start, and a projected gradient-based optimization method.
These methodologies are complemented by a stability-driven early stopping strategy that ensures that the training process returns the best-performing \dISS{} MGU network, if available.
Numerical results showed that the MGU network exhibits superior parameter efficiency and faster inference times compared to the GRU network, while achieving comparable predictive performance. 
This favors the use of the MGU network in embedded predictive control applications, where model stability and resource constraints are of paramount importance.
%
%Comparing the MGU against the GRU (with the same number of hidden units) on the pH neutralization process benchmark, we demonstrated the MGU's superior parameter efficiency and faster inference times compared to the GRU, while achieving comparable predictive performance.
%
Our analysis of stability promotion revealed that, while the combination of loss augmentation and parameters warm-start yielded the best trade-off between prediction accuracy and \dISS{} compliance for the MGU network, loss augmentation alone was insufficient to favor stability consistently for this smaller architecture, even though it aligned with stable GRU network levels of accuracy.
This contrasts with the GRU network, where loss augmentation alone sufficed, likely due to its higher parameter count and over-parameterized nature~\cite[Chapter 8]{goodfellowDeepLearning2016}. 
This finding motivates the design of ad-hoc stability-promoting methodologies for smaller deep learning models, for which the \dISS{} property may not be easily guaranteed. 
Finally, on real-world data from the Silverbox benchmark~\cite{wigrenThreeFreeData2013}, the stable MGU networks significantly outperformed the stable GRU model baseline, demonstrating superior convergence capabilities under stability constraints.
%
%We also found PGD to be unfavorable as it drastically reduced prediction power despite ensuring high \dISS{} compliance.
%
%Ultimately, we have established the MGU as a certifiably stable and highly efficient model for system identification.
\appendix{}
\section{Proofs}
\label{s:proofs}
\paragraph*{Proof of Proposition~\ref{prop:invariant_state_layer}.}
Note that, according to~\eqref{eq:MGU_layer} and for any $l \in \layers$, each component of the hidden state $\boldsymbol{h}_k^{(l)}$ of a single-layer MGU network is updated as:
\begin{tequation}
    \label{eq:MGU_layer_component}
    h_{j,k+1}^{(l)} \! \defeq{} \! \left(1 \! - \! f_{j,k}^{(l)} \right) h_{j,k}^{(l)} \! + \! f_{j,k}^{(l)} \tilde{h}_{j,k}^{(l)}, \quad \forall j \! \in \! \{1, \ldots, n^{(l)}_\hiddenunits{}\}.
\end{tequation}
By definition of the sigmoid and hyperbolic tangent functions (Section~\ref{ss:notation_preliminaries}), $f_{j,k}^{(l)} \in (0,1)$ and $\tilde{h}_{j,k}^{(l)} \in (-1,1)$ regardless of the magnitude of the input $\boldsymbol{\tilde{u}}_k^{(l)}$ in~\eqref{eq:gates_MGU}.
Then, the update equation in~\eqref{eq:MGU_layer_component} is a convex combination.
Since convex combinations preserve membership in convex sets, if $h_{j,k}^{(l)} \in [-1, 1]$, then $h_{j,k+1}^{(l)} \in [-1, 1]$.
By induction from Assumption~\ref{as:initial_state_boundedness}, the hidden state vector $\boldsymbol{h}_{k}^{(l)}$ remains in $\mathcal{H}_{\mathrm{inv}}^{(l)} \defeq{} \left[-1,1\right]^{n_{\hiddenunits}^{(l)}}$ for all $k\in\naturalset$ and any $\boldsymbol{\tilde{u}}_k^{(l)} \in \realset^{n_{\tilde{u}}^{(l)}}$, making $\mathcal{H}_{\mathrm{inv}}^{(l)}$ a forward invariant compact set.
\hfill$\blacksquare$

%%%%%%%%%%%%%%%%%%%%%%%%%%%%%%%%%%%%%%%%%%%%%%%%%%%%%%%%%%%%%%%%%%%%%%%%%%%%%%%%%%%%%%%%%%%%%%%%%%%%%%%%%%

\paragraph*{Proof of Theorem~\ref{th:MGU_layer_ISS}.}
%\label{app:proof_MGU_layer_ISS}
To begin with the proof, we derive an upper bound on the $j$-th component of the forget gate $\boldsymbol{f}_{k}^{(l)}$, $j \in \{1, \ldots, n^{(l)}_\hiddenunits\}$, $l \in \layers$, in~\eqref{eq:gates_MGU} at any $k \in \naturalset$ in a fashion similar to~\cite{terzi_learning_2021,bonassi_stability_2021}, i.e., by leveraging the fact that the sigmoid function is a positive monotonically increasing function:
\begin{talign}
    \label{eq:forget_gate_component_upper_bound}
    \left|f_{j, k}^{(l)}\right| &\leq \infnorm{\boldsymbol{\sigma}\left( W_f^{(l)} {\boldsymbol{\tilde{u}}}_k^{(l)}  +  R_f^{(l)} \boldsymbol{h}_k^{(l)}  +  \boldsymbol{b}_f^{(l)} \right)} \nonumber \\
    &\leq \maxuandh{} \infnorm{\boldsymbol{\sigma}\left( W_f^{(l)} {\boldsymbol{\tilde{u}}}_k^{(l)}  +  R_f^{(l)} \boldsymbol{h}_k^{(l)}  +  \boldsymbol{b}_f^{(l)} \right)} \nonumber \\
    &\leq \sigma \left(\maxuandh{} \infnorm{W_f^{(l)} {\boldsymbol{\tilde{u}}}_k^{(l)}  +  R_f^{(l)} \boldsymbol{h}_k^{(l)}  +  \boldsymbol{b}_f^{(l)}} \right) \nonumber  \\
    &\leq \sigma\left(\infnorm{ \begin{bmatrix}W_{f}^{(l)} & R_{f}^{(l)} & \boldsymbol{b}_{f}^{(l)}\end{bmatrix}}\right) \defeq \bar{\sigma}_{f}^{(l)},
\end{talign}
obtaining the expression in Theorem~\ref{th:MGU_layer_ISS}.
Clearly, it follows that $\infnorm{\boldsymbol{f}_{k}^{(l)}} \leq \bar{\sigma}_{f}^{(l)}$.
Next, owing to the $1$-Lipschitz continuity of the hyperbolic tangent function and using~\eqref{eq:forget_gate_component_upper_bound}, we obtain an upper bound on the candidate hidden state $\boldsymbol{\tilde{h}}_{k}^{(l)}$ in~\eqref{eq:gates_MGU} as
\begin{talign}
    \label{eq:candidate_hidden_state_upper_bound}
    \infnorm{\boldsymbol{\tilde{h}}_{k}^{(l)}} &\leq \infnorm{W_{\tilde{h}}^{(l)} {\boldsymbol{\tilde{u}}}_k^{(l)}  +  R_{\tilde{h}}^{(l)} \left( \boldsymbol{f}_k^{(l)}  \circ  \boldsymbol{h}_k^{(l)} \right)  +  \boldsymbol{b}_{\tilde{h}}^{(l)}} \nonumber\\
    &\leq \infnorm{W_{\tilde{h}}^{(l)}} \infnorm{\boldsymbol{\tilde{u}}_k^{(l)}} + \infnorm{R_{\tilde{h}}^{(l)}} \infnorm{\boldsymbol{f}_k^{(l)}} \infnorm{\boldsymbol{h}_k^{(l)}} + \infnorm{\boldsymbol{b}_{\tilde{h}}^{(l)}} \nonumber \\
    & \leq \infnorm{W_{\tilde{h}}^{(l)}} \infnorm{\boldsymbol{\tilde{u}}_k^{(l)}} + \bar{\sigma}_{f}^{(l)} \infnorm{R_{\tilde{h}}^{(l)}} \infnorm{\boldsymbol{h}_k^{(l)}} + \infnorm{\boldsymbol{b}_{\tilde{h}}^{(l)}}.
\end{talign}
Next, leveraging~\eqref{eq:forget_gate_component_upper_bound} and~\eqref{eq:candidate_hidden_state_upper_bound}, we bound the $j$-th component of the hidden state in~\eqref{eq:MGU_layer_component} as\footnote{Note that $f_{j,k}^{(l)} = \left|f_{j,k}^{(l)}\right|$ and $1 - f_{j,k}^{(l)} = \left|1 - f_{j,k}^{(l)}\right|$ since the range of the sigmoid function is $(0, 1)$.}:
\begin{talign}
    \label{eq:layer_ISS_intermediate_step_1}
    \left|h_{j,k+1}^{(l)}\right| &\leq \left(1  -  f_{j,k}^{(l)} \right) \left|h_{j,k}^{(l)}\right| +  f_{j,k}^{(l)} \left|\tilde{h}_{j,k}^{(l)}\right| \nonumber \\
    &\leq \left(1  -  f_{j,k}^{(l)} \right) \infnorm{\boldsymbol{h}_k^{(l)}} +  f_{j,k}^{(l)} \infnorm{\boldsymbol{\tilde{h}}_{k}^{(l)}} \nonumber \\
    & \leq \left(1  -  f_{j,k}^{(l)} + f_{j,k}^{(l)}  \bar{\sigma}_{f}^{(l)} \infnorm{R_{\tilde{h}}^{(l)}}\right) \infnorm{\boldsymbol{h}_k^{(l)}} + \\
    & \quad + f_{j,k}^{(l)} \left(\infnorm{W_{\tilde{h}}^{(l)}} \infnorm{\boldsymbol{\tilde{u}}_k^{(l)}} + \infnorm{\boldsymbol{b}_{\tilde{h}}^{(l)}}\right). \nonumber
\end{talign}
Now, consider the condition in~\eqref{eq:MGU_layer_ISS_condition} and note that $f_{j,k}^{(l)} \in (0, 1)$ due to the range of the sigmoid function.
Then, the term $1  -  f_{j,k}^{(l)} + f_{j,k}^{(l)} \bar{\sigma}_{f}^{(l)} \infnorm{R_{\tilde{h}}^{(l)}}$ in~\eqref{eq:layer_ISS_intermediate_step_1} is a convex combination between $1$ and $\bar{\sigma}_{f}^{(l)} \infnorm{R_{\tilde{h}}^{(l)}}$, which must necessarily be in $(0, 1)$.
Hence, for any $j \in \{1, \ldots, n^{(l)}_\hiddenunits\}$, there must exist a $\lambda_j^{(l)} \in (0, 1)$ such that $1  -  f_{j,k}^{(l)} + f_{j,k}^{(l)} \bar{\sigma}_{f}^{(l)} \infnorm{R_{\tilde{h}}^{(l)}} \leq \lambda_j^{(l)}$.
Defining $\lambda^{(l)} \defeq{} \max_{j \in \{1, \ldots, n^{(l)}_\hiddenunits\}} \lambda_j^{(l)}$ and returning to~\eqref{eq:layer_ISS_intermediate_step_1}, we have:
\begin{talign*}
    \left|h_{j,k+1}^{(l)}\right| &\leq \lambda^{(l)} \infnorm{\boldsymbol{h}_k^{(l)}} + f_{j,k}^{(l)} \left(\infnorm{W_{\tilde{h}}^{(l)}} \infnorm{\boldsymbol{\tilde{u}}_k^{(l)}} + \infnorm{\boldsymbol{b}_{\tilde{h}}^{(l)}}\right) \\
    &\leq \lambda^{(l)} \infnorm{\boldsymbol{h}_k^{(l)}} + \bar{\sigma}_{f}^{(l)} \left(\infnorm{W_{\tilde{h}}^{(l)}} \infnorm{\boldsymbol{\tilde{u}}_k^{(l)}} + \infnorm{\boldsymbol{b}_{\tilde{h}}^{(l)}}\right).
\end{talign*}
Taking the maximum w.r.t. $j \in \{1, \ldots, n^{(l)}_\hiddenunits\}$ on both sides and shifting the time index backwards by one step, we get:
\begin{tequation}
    \label{eq:hidden_state_upper_bound}
    \infnorm{\boldsymbol{h}_{k}^{(l)}} \leq \alpha^{(l)} \infnorm{\boldsymbol{h}_{k-1}^{(l)}} + \beta_u^{(l)} \infnorm{\boldsymbol{\tilde{u}}_{k-1}^{(l)}} +  \beta_b^{(l)}  \infnorm{\boldsymbol{b}_{\tilde{h}}^{(l)}},
\end{tequation}
where $\alpha^{(l)} = \lambda^{(l)}$, $\beta_u^{(l)} = \bar{\sigma}_{f}^{(l)} \infnorm{W_{\tilde{h}}^{(l)}}$, and $\beta_b^{(l)} = \bar{\sigma}_{f}^{(l)}$.
Next, we propagate~\eqref{eq:hidden_state_upper_bound} backwards to the initial state $\boldsymbol{h}_{0}^{(l)} \in \mathcal{H}_{\mathrm{inv}}^{(l)}$:
\begin{talign}
    \label{eq:hidden_state_upper_bound_propagated}
    \infnorm{\boldsymbol{h}_{k}^{(l)}} \! &\leq \! \alpha^{(l)^k} \! \infnorm{\boldsymbol{h}_0^{(l)}} \! + \! \sum_{z=0}^{k-1} \alpha^{(l)^{k - 1 - z}} \!\!\left(\beta_u^{(l)} \! \infnorm{\boldsymbol{\tilde{u}}_{z}^{(l)}} \! + \! \beta_b^{(l)} \! \infnorm{\boldsymbol{b}_{\tilde{h}}^{(l)}}\right) \nonumber \\
    &\leq \! \alpha^{(l)^k} \! \infnorm{\boldsymbol{h}_0^{(l)}} \! + \! \beta_u^{(l)} \max_{0 \leq z < k} \infnorm{\boldsymbol{\tilde{u}}_{z}^{(l)}} \sum_{z=0}^{k-1} \alpha^{(l)^{k - 1 - z}} \! + \nonumber \\
    &\quad + \! \beta_b^{(l)} \infnorm{\boldsymbol{b}_{\tilde{h}}^{(l)}} \sum_{z=0}^{k-1} \alpha^{(l)^{k - 1 - z}}.
\end{talign}
Given that $\lambda^{(l)} \in (0, 1)$, the geometric series $\sum_{z=0}^{k-1} \alpha^{(l)^{k - 1 - z}}$ can be bounded as
\begin{tequation}
    \label{eq:convergent_geometric_series_proof}
    \sum_{z=0}^{k-1} \alpha^{(l)^{k - 1 - z}} = \frac{1 - \lambda^{(l)^k}}{1 - \lambda^{(l)}} \leq \frac{1}{1 - \lambda^{(l)}}.
\end{tequation}
Finally, still owing to the fact that $\lambda^{(l)} \in (0, 1)$ and leveraging~\eqref{eq:convergent_geometric_series_proof}, we can immediately deduce the $\mathcal{KL}$ function for the ISS condition in~\eqref{eq:iss} from~\eqref{eq:hidden_state_upper_bound_propagated}, i.e.,
\begin{tequation}
    \label{eq:KL_ISS_layer-wise}
    \psi\left(\infnorm{\boldsymbol{h}_{0}^{(l)}}, k\right) \defeq{} \lambda^{(l)^k} \infnorm{\boldsymbol{h}_0^{(l)}},
\end{tequation}
and the two $\mathcal{K}_{\infty}$ functions, i.e.,
\begin{talign}
    \label{eq:Kinfty_ISS_layer-wise_input}
    \gamma_u \left(\max_{0 \leq z < k} \infnorm{\boldsymbol{\tilde{u}}_{z}^{(l)}} \right) &\defeq{} \frac{\bar{\sigma}_{f}^{(l)} \infnorm{W_{\tilde{h}}^{(l)}}}{1 - \lambda^{(l)}} \max_{0 \leq z < k} \infnorm{\boldsymbol{\tilde{u}}_{z}^{(l)}}, \\
    \label{eq:Kinfty_ISS_layer-wise_bias}
    \gamma_b \left(\infnorm{\boldsymbol{b}_{\tilde{h}}^{(l)}} \right)  &\defeq{} \frac{\bar{\sigma}_{f}^{(l)}}{1 - \lambda^{(l)}} \infnorm{\boldsymbol{b}_{\tilde{h}}^{(l)}}.
\end{talign}
We can conclude that, under Condition~\eqref{eq:MGU_layer_ISS_condition}, the $l$-th, $l \in \layers$, MGU layer is ISS in $\mathcal{H}_{\mathrm{inv}}^{(l)}$ with respect to $\tilde{\mathcal{U}}^{(l)}$ due to the existence of the functions~\eqref{eq:KL_ISS_layer-wise}~-~\eqref{eq:Kinfty_ISS_layer-wise_bias}.
\hfill$\blacksquare$

%%%%%%%%%%%%%%%%%%%%%%%%%%%%%%%%%%%%%%%%%%%%%%%%%%%%%%%%%%%%%%%%%%%%%%%%%%%%%%%%%%%%%%%%%%%%%%%%%%%%%%%%%%

\paragraph*{Proof of Theorem~\ref{th:MGU_layer_dISS}}
%\label{app:proof_MGU_layer_dISS}
Let us consider two state trajectories $\boldsymbol{h}_k^{(l), \mathrm{a}}$ and $\boldsymbol{h}_k^{(l), \mathrm{b}}$ for the $l$-th, $l \in \layers$, MGU layer in~\eqref{eq:MGU_layer} driven by the input sequences $\{\boldsymbol{\tilde{u}}^{(l), \mathrm{a}}_{z}\in\tilde{\mathcal{U}}^{(l)}\}_{z=0}^{k-1}$ and $\{\boldsymbol{\tilde{u}}^{(l), \mathrm{b}}_{z}\in\tilde{\mathcal{U}}^{(l)}\}_{z=0}^{k-1}$, respectively, and started from the initial hidden states $\boldsymbol{h}_0^{(l), \mathrm{a}} \in \mathcal{H}_{\mathrm{inv}}^{(l)}$ and $\boldsymbol{h}_0^{(l), \mathrm{b}} \in \mathcal{H}_{\mathrm{inv}}^{(l)}$, respectively.
Throughout this proof, we will use the superscript $\Delta$ to refer to the incremental dynamics, e.g., $\boldsymbol{h}_k^{(l), \Delta} = \boldsymbol{h}_k^{(l), \mathrm{a}} - \boldsymbol{h}_k^{(l), \mathrm{b}}$.
The incremental dynamics for the $j$-th component, $j \in \{1, \ldots, n^{(l)}_\hiddenunits\}$, of the $l$-th hidden state in~\eqref{eq:MGU_layer_component} amount to:
\begin{talign*}
    h_{j, k+1}^{(l), \Delta} &= \left(1-f_{j,k}^{\left(l\right),\mathrm{a}}\right)h_{j,k}^{\left(l\right),\mathrm{a}}-\left(1-f_{j,k}^{\left(l\right),\mathrm{b}}\right)h_{j,k}^{\left(l\right),\mathrm{b}} + \\
    &\quad + f_{j,k}^{\left(l\right),\mathrm{a}}\tilde{h}_{j,k}^{\left(l\right),\mathrm{a}}-f_{j,k}^{\left(l\right),\mathrm{b}}\tilde{h}_{j,k}^{\left(l\right),\mathrm{b}}.
\end{talign*}
Adding and subtracting $\left(1-f_{j,k}^{\left(l\right),\mathrm{a}}\right)h_{j,k}^{\left(l\right),\mathrm{b}}$ and $f_{j,k}^{\left(l\right),\mathrm{a}}\tilde{h}_{j,k}^{\left(l\right),\mathrm{b}}$ to the right side of the previous expression, we get:
\begin{tequation}
    \label{eq:layer_dISS_intermediate_step_1}
    h_{j,k+1}^{\left(l\right),\Delta} \! = \! \left(1 \! - \! f_{j,k}^{\left(l\right),\mathrm{a}}\right) \! h_{j,k}^{\left(l\right),\Delta} \! + \! f_{j,k}^{\left(l\right),\mathrm{a}}\tilde{h}_{j,k}^{\left(l\right),\Delta} \! + \! f_{j,k}^{\left(l\right),\Delta} \!\! \left(\tilde{h}_{j,k}^{\left(l\right),\mathrm{b}} \!\! - \! h_{j,k}^{\left(l\right),\mathrm{b}}\right)\!.
\end{tequation}
Now, note that $\left|f_{j,k}^{\left(l\right),\mathrm{a}}\right| \leq \bar{\sigma}_{f}^{(l)}$ due to~\eqref{eq:forget_gate_component_upper_bound}, and, following similar mathematical steps as in~\eqref{eq:forget_gate_component_upper_bound}, we can obtain the upper bound $\left|1 - f_{j,k}^{\left(l\right),\mathrm{a}}\right| \leq \bar{\sigma}_{f}^{(l)}$.
Then, the incremental component-wise dynamics of the hidden state in~\eqref{eq:layer_dISS_intermediate_step_1} are bounded as follows:
\begin{talign}
    \label{eq:layer_dISS_intermediate_step_2}
   &\left| h_{j,k+1}^{\left(l\right),\Delta} \right| \! \leq \! \left|1 \! - \! f_{j,k}^{\left(l\right),\mathrm{a}}\right| \left| h_{j,k}^{\left(l\right),\Delta} \right| \! + \! \left| f_{j,k}^{\left(l\right),\mathrm{a}} \right| \left| \tilde{h}_{j,k}^{\left(l\right),\Delta} \right| \! + \nonumber \\
    & \qquad \qquad + \! \left| f_{j,k}^{\left(l\right),\Delta} \right| \left| \tilde{h}_{j,k}^{\left(l\right),\mathrm{b}} \! - \! h_{j,k}^{\left(l\right),\mathrm{b}} \right| \\
    %&\leq \bar{\sigma}_{f}^{(l)} \left| h_{j,k}^{\left(l\right),\Delta} \right|  \! + \! \bar{\sigma}_{f}^{(l)} \left| \tilde{h}_{j,k}^{\left(l\right),\Delta} \right| \! + \! \left| f_{j,k}^{\left(l\right),\Delta} \right| \left| \tilde{h}_{j,k}^{\left(l\right),\mathrm{b}} \! - \! h_{j,k}^{\left(l\right),\mathrm{b}} \right|
    &\quad \leq \bar{\sigma}_{f}^{(l)} \! \left| h_{j,k}^{\left(l\right),\Delta} \right|  \! + \! \bar{\sigma}_{f}^{(l)} \! \left| \tilde{h}_{j,k}^{\left(l\right),\Delta} \right| \! + \! \left| f_{j,k}^{\left(l\right),\Delta} \right| \left( \left| \tilde{h}_{j,k}^{\left(l\right),\mathrm{b}} \right| \! + \! \left| h_{j,k}^{\left(l\right),\mathrm{b}} \right| \right). \nonumber 
\end{talign}
The incremental dynamics of the forget gate $\left| f_{j,k}^{\left(l\right),\Delta} \right| \leq \infnorm{\boldsymbol{f}_k^{(l), \Delta}}$ can be upper bounded by leveraging the $\frac{1}{4}$-Lipschitz continuity of the sigmoid function, leading to:
\begin{talign}
    \label{eq:forget_gate_incremental_upper_bound}
    %\infnorm{\boldsymbol{f}_k^{(l), \Delta}} \! \leq \!\frac{1}{4} \infnorm{W_f^{(l)} {\boldsymbol{\tilde{u}}}_k^{(l), \mathrm{a}}  \! + \! R_f^{(l)} \boldsymbol{h}_k^{(l), \mathrm{a}}  \! + \! \boldsymbol{b}_f^{(l)} \! - \! W_f^{(l)} {\boldsymbol{\tilde{u}}}_k^{(l), \mathrm{b}} \! - \! R_f^{(l)} \boldsymbol{h}_k^{(l), \mathrm{b}} \! + \! \boldsymbol{b}_f^{(l)}} \nonumber
    \!\!\infnorm{\boldsymbol{f}_k^{(l), \Delta}} \! \leq \! \frac{1}{4} \! \infnorm{W_f^{(l)}} \! \infnorm{\boldsymbol{\tilde{u}}_k^{(l), \Delta}} \! + \! \frac{1}{4} \! \infnorm{R_f^{(l)}} \! \infnorm{\boldsymbol{h}_k^{(l), \Delta}}\!.
\end{talign}
Similarly, the incremental dynamics of the candidate hidden state $\left| \tilde{h}_{j,k}^{\left(l\right),\Delta} \right| \leq \infnorm{\boldsymbol{\tilde{h}}_k^{(l), \Delta}}$ can be bounded by exploiting the $1$-Lipschitz continuity of the hyperbolic tangent function:
\begin{talign*}
    \infnorm{\boldsymbol{\tilde{h}}_k^{(l), \Delta}} \! &\leq \! \infnorm{W_{\tilde{h}}^{(l)}} \! \infnorm{\boldsymbol{\tilde{u}}_k^{(l), \Delta}} \! + \\
    &\quad + \! \infnorm{R_{\tilde{h}}^{(l)}} \!\infnorm{\boldsymbol{f}_k^{(l), \mathrm{a}} \!\! \circ \! \boldsymbol{h}_k^{(l), \mathrm{a}} \!\! - \! \boldsymbol{f}_k^{(l), \mathrm{b}} \!\! \circ \! \boldsymbol{h}_k^{(l), \mathrm{b}}}\!.
\end{talign*}
By adding and subtracting $\boldsymbol{f}_{k}^{(l),\mathrm{a}}\circ\boldsymbol{h}_{k}^{(l),\mathrm{b}}$ within the right-most norm of the previous expression and substituting~\eqref{eq:forget_gate_incremental_upper_bound}, taking also into account that $\infnorm{\boldsymbol{h}_k^{(l), \mathrm{b}}} \leq 1$ since $\boldsymbol{h}_k^{(l), \mathrm{b}} \in \mathcal{H}_{\mathrm{inv}}^{(l)}$, we get:
\begin{talign}
    \label{eq:incremental_candidate_hidden_state_bound}
    \infnorm{\boldsymbol{\tilde{h}}_k^{(l), \Delta}} \! &\leq \! \infnorm{W_{\tilde{h}}^{(l)}} \! \infnorm{\boldsymbol{\tilde{u}}_k^{(l), \Delta}} \! + \! \bar{\sigma}_{f}^{(l)} \! \infnorm{R_{\tilde{h}}^{(l)}} \infnorm{\boldsymbol{h}_k^{(l), \Delta}} + \nonumber\\
    &\quad + \! \infnorm{R_{\tilde{h}}^{(l)}} \infnorm{\boldsymbol{f}_k^{(l), \Delta}} \\
    &\leq \! \left(\infnorm{W_{\tilde{h}}^{(l)}} \! + \! \frac{1}{4} \! \infnorm{R_{\tilde{h}}^{(l)}} \! \infnorm{W_f^{(l)}}\right) \! \infnorm{\boldsymbol{\tilde{u}}_k^{(l), \Delta}} \! + \nonumber \\
    & \quad + \! \left(\bar{\sigma}_{f}^{(l)} \! \infnorm{R_{\tilde{h}}^{(l)}} \! + \! \frac{1}{4} \! \infnorm{R_{\tilde{h}}^{(l)}} \! \infnorm{R_f^{(l)}}\right) \! \infnorm{\boldsymbol{h}_k^{(l), \Delta}}. \nonumber
\end{talign}
Next, we bound the term $\left| \tilde{h}_{j,k}^{\left(l\right),\mathrm{b}} \right|$ in~\eqref{eq:layer_dISS_intermediate_step_2}.
For this purpose, recall that $\phi(|x|) = |\phi(x)|$ for any $x \in \mathbb{R}$ and $\phi(|x|)$ is a monotonically increasing function.
Then:
\begin{talign}
    \label{eq:candidate_hidden_state_bound} 
    \left| \tilde{h}_{j,k}^{\left(l\right),\mathrm{b}} \right| \! &\leq \!\infnorm{\boldsymbol{\phi} \! \left( W_{\tilde{h}}^{(l)} {\boldsymbol{\tilde{u}}}_k^{(l), \mathrm{b}}  \!\! + \! R_{\tilde{h}}^{(l)} \! \left( \boldsymbol{f}_k^{(l), \mathrm{b}} \!\! \circ \! \boldsymbol{h}_k^{(l), \mathrm{b}} \right)  \!\! + \! \boldsymbol{b}_{\tilde{h}}^{(l)} \right)} \nonumber \\
    &= \! \phi\left(\infnorm{W_{\tilde{h}}^{(l)} {\boldsymbol{\tilde{u}}}_k^{(l), \mathrm{b}}  \!\! + \! R_{\tilde{h}}^{(l)} \! \left( \boldsymbol{f}_k^{(l), \mathrm{b}} \!\! \circ \! \boldsymbol{h}_k^{(l), \mathrm{b}} \right)  \!\! + \! \boldsymbol{b}_{\tilde{h}}^{(l)}}\right) \nonumber \\
    &\leq \!\!\! \maxuandhstackedb \!\!\!\! \phi\left(\infnorm{W_{\tilde{h}}^{(l)} {\boldsymbol{\tilde{u}}}_k^{(l), \mathrm{b}}  \!\! + \! R_{\tilde{h}}^{(l)} \!\! \left( \boldsymbol{f}_k^{(l), \mathrm{b}} \!\! \circ \! \boldsymbol{h}_k^{(l), \mathrm{b}} \right)  \!\! + \! \boldsymbol{b}_{\tilde{h}}^{(l)}}\right) \nonumber \\
    &= \! \phi \! \left(\!\maxuandhstackedb \!\!\! \infnorm{W_{\tilde{h}}^{(l)} {\boldsymbol{\tilde{u}}}_k^{(l), \mathrm{b}}  \!\! + \! R_{\tilde{h}}^{(l)} \!\! \left( \boldsymbol{f}_k^{(l), \mathrm{b}} \!\! \circ \! \boldsymbol{h}_k^{(l), \mathrm{b}} \right)  \!\! + \! \boldsymbol{b}_{\tilde{h}}^{(l)}} \!\!\right) \nonumber \\
    &= \phi \left(\infnorm{W_{\tilde{h}}^{(l)} \boldsymbol{1}_{n_{\tilde{u}}^{(l)}} \! + \! R_{\tilde{h}}^{(l)} \boldsymbol{f}_k^{(l), \mathrm{b}} \! + \! \boldsymbol{b}_{\tilde{h}}^{(l)}}\right) \nonumber \\
    &\leq \phi\left(\infnorm{ \begin{bmatrix}W_{\tilde{h}}^{(l)} & R_{\tilde{h}}^{(l)} & \boldsymbol{b}_{\tilde{h}}^{(l)}\end{bmatrix}} \right) \defeq{} \bar\phi_{\tilde{h}}^{(l)},
\end{talign}
obtaining the expression in Theorem~\ref{th:MGU_layer_dISS}.
It immediately follows that $\infnorm{\boldsymbol{\tilde{h}}_k^{(l), \mathrm{b}}} \leq \bar\phi_{\tilde{h}}^{(l)}$.
The last term in~\eqref{eq:layer_dISS_intermediate_step_2} that needs to be upper bounded is $\left| h_{j,k}^{\left(l\right),\mathrm{b}} \right|$.
Due to the forward invariance of $\mathcal{H}_{\mathrm{inv}}^{(l)}$ (Proposition~\ref{prop:invariant_state_layer}) and Assumption~\ref{as:initial_state_boundedness}, we can immediately conclude that $\left| h_{j,k}^{\left(l\right),\mathrm{b}} \right| \leq 1$.
Now, going back to~\eqref{eq:layer_dISS_intermediate_step_2} and leveraging~\eqref{eq:forget_gate_incremental_upper_bound}~-~\eqref{eq:candidate_hidden_state_bound}, we get:
\begin{talign}
    \label{eq:layer_dISS_intermediate_step_3}
   \left| h_{j,k+1}^{\left(l\right),\Delta} \right| \! &\leq \! \bar{\sigma}_{f}^{(l)} \!  \infnorm{\boldsymbol{h}_k^{(l), \Delta}} \! + \! \bar{\sigma}_{f}^{(l)} \! \infnorm{\boldsymbol{\tilde{h}}_k^{(l), \Delta}} \! + \! \infnorm{\boldsymbol{f}_k^{(l), \Delta}} \!\! \left( \bar\phi_{\tilde{h}}^{(l)}\! + \! 1\right) \nonumber  \\
   &\leq \alpha_{\delta}^{(l)} \infnorm{\boldsymbol{h}_k^{(l), \Delta}} + \beta_{\delta u}^{(l)} \infnorm{\boldsymbol{\tilde{u}}_k^{(l), \Delta}},
\end{talign}
where
\begin{talign*}
    \alpha_{\delta}^{(l)} \! &= \! \bar{\sigma}_{f}^{(l)} \! + \! \bar{\sigma}_{f}^{(l)^{2}} \! \infnorm{ R_{\tilde{h}}^{(l)}} \! + \! \frac{1}{4} \infnorm{R_{f}^{(l)}} \! \left(\bar{\sigma}_{f}^{(l)} \infnorm{R_{\tilde{h}}^{(l)}} \! + \! \bar\phi_{\tilde{h}}^{(l)} \! + \! 1 \! \right), \\
    \beta_{\delta u}^{(l)} \! &= \! \bar{\sigma}_{f}^{(l)} \! \infnorm{W_{\tilde{h}}^{(l)}} \! + \! \frac{1}{4} \! \infnorm{W_{f}^{(l)}} \! \left(\bar{\sigma}_{f}^{(l)} \! \infnorm{ R_{\tilde{h}}^{(l)}} \! + \! \bar\phi_{\tilde{h}}^{(l)}\! + \! 1\right).
\end{talign*}
Note that $\alpha_{\delta}^{(l)}$ corresponds to the left-hand side of the condition in~\eqref{eq:MGU_layer_dISS_condition}.
Next, we take the maximum w.r.t. $j \in \{1, \ldots, n^{(l)}_\hiddenunits\}$ on both sides of~\eqref{eq:layer_dISS_intermediate_step_3} and shift the time index backwards by one step, obtaining:
\begin{tequation}
    \label{eq:incremental_hidden_state_upper_bound}
    \infnorm{\boldsymbol{h}_{k}^{(l), \Delta}} \leq \alpha_{\delta}^{(l)} \infnorm{\boldsymbol{h}_{k-1}^{(l), \Delta}} + \beta_{\delta u}^{(l)} \infnorm{\boldsymbol{\tilde{u}}_{k-1}^{(l)}}.
\end{tequation}
Finally, we propagate~\eqref{eq:incremental_hidden_state_upper_bound} backwards to the initial state difference $\boldsymbol{h}_{0}^{(l), \Delta}$ in a fashion similar to~\eqref{eq:hidden_state_upper_bound_propagated} to get:
\begin{tequation}
    \label{eq:incrementa_hidden_state_upper_bound_propagated}
    \infnorm{\boldsymbol{h}_{k}^{(l), \Delta}} \leq \alpha_{\delta}^{(l)^k}  \infnorm{\boldsymbol{h}_0^{(l), \Delta}} +  \beta_{\delta u}^{(l)} \max_{0 \leq z < k} \infnorm{\boldsymbol{\tilde{u}}_{z}^{(l), \Delta}} \sum_{z=0}^{k-1} \alpha_{\delta}^{(l)^{k - 1 - z}} \!\!\!.
\end{tequation}
Under Condition~\eqref{eq:MGU_layer_dISS_condition}, the right-most summation in the previous expression can be upper bounded as in~\eqref{eq:convergent_geometric_series_proof}, and we can consequently find the $\mathcal{KL}$ and $\mathcal{K}_{\infty}$ functions for the \dISS{} condition in~\eqref{eq:diss}:
\begin{talign}
    \label{eq:KL_dISS_layer-wise}
    &\!\!\psi_{\delta} \! \left(\!\infnorm{\boldsymbol{h}_{0}^{(l),\mathrm{a}} \!\! - \! \boldsymbol{h}_{0}^{(l),\mathrm{b}}} \! , k\!\right) \! \defeq{} \! \alpha_{\delta}^{\left(l\right)^{k}} \!\infnorm{\boldsymbol{h}_{0}^{(l),\mathrm{a}} \!\! - \! \boldsymbol{h}_{0}^{(l),\mathrm{b}}} \!, \\
    \label{eq:Kinfty_dISS_layer-wise_input}
    &\!\!\gamma_{\delta u} \! \left(\!\max_{0 \leq z <k } \! \infnorm{\boldsymbol{\tilde{u}}_{z}^{\left(l\right),\mathrm{a}} \!\! - \! \boldsymbol{\tilde{u}}_{z}^{\left(l\right),\mathrm{b}} }\!\right) \! \defeq \! \frac{\beta_{\delta u}^{\left(l\right)}}{1 \! - \! \alpha_{\delta}^{\left(l\right)}} \! \max_{0\leq z<k} \! \infnorm{\boldsymbol{\tilde{u}}_{z}^{\left(l\right),\mathrm{a}} \!\! - \!\boldsymbol{\tilde{u}}_{z}^{\left(l\right),\mathrm{b}}} \!.
\end{talign}
We can conclude that, under Condition~\eqref{eq:MGU_layer_dISS_condition}, the $l$-th, $l \in \layers$, MGU layer is \dISS{} in $\mathcal{H}_{\mathrm{inv}}^{(l)}$ with respect to $\tilde{\mathcal{U}}^{(l)}$ due to the existence of the functions~\eqref{eq:KL_dISS_layer-wise}~and~\eqref{eq:Kinfty_dISS_layer-wise_input}.
\hfill$\blacksquare$

%%%%%%%%%%%%%%%%%%%%%%%%%%%%%%%%%%%%%%%%%%%%%%%%%%%%%%%%%%%%%%%%%%%%%%%%%%%%%%%%%%%%%%%%%%%%%%%%%%%%%%%%%%

\paragraph*{Proof of Proposition~\ref{pr:MGU_layer_dISS_imply_ISS}}
%\label{app:proof_MGU_layer_dISS_imply_ISS}
%
Consider the \dISS{} condition in~\eqref{eq:MGU_layer_dISS_condition} for the $l$-th MGU layer, $l \in \layers$, and let $\varsigma^{\left(l\right)} \defeq{} \bar{\sigma}_{f}^{\left(l\right)} \infnorm{R_{\tilde{h}}^{(l)}}$, amounting to the left-hand side of Condition~\eqref{eq:MGU_layer_ISS_condition}.
We can re-write Condition~\eqref{eq:MGU_layer_dISS_condition} as follows:
\begin{tequation*}
    c_1^{(l)} \!\! \defeq{} \! \bar{\sigma}_{f}^{\left(l\right)} + \bar{\sigma}_{f}^{\left(l\right)}\varsigma^{\left(l\right)} + \frac{1}{4}\! \infnorm{R_{f}^{\left(l\right)}} \varsigma^{\left(l\right)} + \frac{1}{4} \! \infnorm{R_{f}^{\left(l\right)}} \bar{\phi}_{\tilde{h}}^{\left(l\right)} + \frac{1}{4} \! \infnorm{R_{f}^{\left(l\right)}} \! < \! 1.
\end{tequation*}
Given that $\frac{1}{4} \! \infnorm{R_{f}^{\left(l\right)}} \bar{\phi}_{\tilde{h}}^{\left(l\right)} \in \realset_{\geq 0}$ since $\infnorm{R_{f}^{\left(l\right)}} \in \realset_{\geq 0}$ and $\bar{\phi}_{\tilde{h}}^{\left(l\right)} \in [0, 1)$, we can derive a lower bound for $c_1^{(l)}$ by omitting the just-mentioned term:
\begin{tequation*}
    c_2^{(l)} \! \defeq{} \! \bar{\sigma}_{f}^{\left(l\right)} \! + \! \bar{\sigma}_{f}^{\left(l\right)}\varsigma^{\left(l\right)} \! + \! \frac{1}{4}\! \infnorm{R_{f}^{\left(l\right)}} \varsigma^{\left(l\right)} \! + \! \frac{1}{4} \! \infnorm{R_{f}^{\left(l\right)}} \! \leq \! c_1^{(l)} \! < \! 1.
\end{tequation*}
Then, the terms that compose $c_2^{(l)}$ can be re-organized to get:
\begin{tequation*}
    c_2^{(l)} = \left(\bar{\sigma}_{f}^{\left(l\right)} + \frac{1}{4} \! \infnorm{R_{f}^{\left(l\right)}}\right) \left( 1 + \varsigma^{\left(l\right)}\right).
\end{tequation*}
Now, note that $\bar{\sigma}_{f}^{\left(l\right)}$ in~\eqref{eq:forget_gate_component_upper_bound} is lower bounded by $\sigma(0) = \frac{1}{2}$ while $\infnorm{R_{f}^{\left(l\right)}} \in \realset_{\geq 0}$.
Hence, we can derive the following lower bound for $c_2^{(l)}$:
\begin{tequation*}
    c_3^{(l)} \defeq{} \frac{1}{2} \left(1 + \varsigma^{\left(l\right)}\right) \leq c_2^{(l)} \leq c_1^{(l)} < 1.
\end{tequation*}
In turn, this means that $\varsigma^{\left(l\right)} < 1$, which is exactly the ISS condition in~\eqref{eq:MGU_layer_ISS_condition}.
Thus, satisfaction of the \dISS{} condition in~\eqref{eq:MGU_layer_dISS_condition} analytically necessitates the ISS condition in~\eqref{eq:MGU_layer_ISS_condition}.
\hfill$\blacksquare$

%%%%%%%%%%%%%%%%%%%%%%%%%%%%%%%%%%%%%%%%%%%%%%%%%%%%%%%%%%%%%%%%%%%%%%%%%%%%%%%%%%%%%%%%%%%%%%%%%%%%%%%%%%

\paragraph*{Proof of Theorem~\ref{th:MGU_network_ISS}}
Theorem~\ref{th:MGU_layer_ISS} establishes that satisfying Condition~\eqref{eq:MGU_layer_ISS_condition} guarantees the ISS property for any MGU layer $l\in\layers$.
The MGU network is a cascade interconnection of these layers, where $\boldsymbol{h}_{k+1}^{(l-1)}$ acts as the input to layer $l \in \layers \setminus \{1\}$ (see~\eqref{eq:input_for_each_layer}).
As proven in~\cite[Theorem 2]{jiang_input--state_2001}, a cascade of ISS discrete-time subsystems is itself ISS.
Thus, the MGU network is ISS.
\hfill$\blacksquare$

%%%%%%%%%%%%%%%%%%%%%%%%%%%%%%%%%%%%%%%%%%%%%%%%%%%%%%%%%%%%%%%%%%%%%%%%%%%%%%%%%%%%%%%%%%%%%%%%%%%%%%%%%%

\paragraph*{Proof of Theorem~\ref{th:MGU_network_dISS}}
%\label{app:proof_MGU_network_dISS}
Let us consider two state trajectories $\boldsymbol{h}_k^{\mathrm{a}} = \left[\boldsymbol{h}_k^{(1), \mathrm{a}^\top}, \ldots, \boldsymbol{h}_k^{(L), \mathrm{a}^\top}\right]^\top$ and $\boldsymbol{h}_k^{\mathrm{b}} = \left[\boldsymbol{h}_k^{(1), \mathrm{b}^\top}, \ldots, \boldsymbol{h}_k^{(L), \mathrm{b}^\top}\right]^\top$ for the MGU network in~\eqref{eq:MGU_network} driven by the input sequences $\{\boldsymbol{u}_z^{\mathrm{a}} \in \mathcal{U}\}_{z=0}^{k-1}$ and $\{\boldsymbol{u}_z^{\mathrm{b}} \in \mathcal{U}\}_{z=0}^{k-1}$, respectively, and started from the initial hidden states $\boldsymbol{h}_0^{\mathrm{a}} \in \mathcal{H}_{\mathrm{inv}}$ and $\boldsymbol{h}_0^{\mathrm{b}} \in \mathcal{H}_{\mathrm{inv}}$, respectively.
Similarly to the proof of Theorem~\ref{th:MGU_layer_dISS}, we will use the superscript $\Delta$ to refer to the incremental dynamics, e.g., $\boldsymbol{h}_k^{\Delta} = \boldsymbol{h}_k^{\mathrm{a}} - \boldsymbol{h}_k^{\mathrm{b}}$ and $\boldsymbol{u}_{k-1}^{\Delta} = \boldsymbol{u}_{k-1}^{\mathrm{a}} - \boldsymbol{u}_{k-1}^{\mathrm{b}}$.
The starting point for this proof is~\eqref{eq:incremental_hidden_state_upper_bound}, which can be unwrapped for each layer $l \in \layers$ by leveraging~\eqref{eq:input_for_each_layer}:
\begin{talign}
    \label{eq:incremental_hidden_state_upper_bound_all}
    \infnorm{\boldsymbol{h}_{k}^{(1), \Delta}} &\leq \alpha_{\delta}^{(1)} \infnorm{\boldsymbol{h}_{k-1}^{(1), \Delta}} + \beta_{\delta u}^{(1)} \infnorm{\boldsymbol{u}_{k-1}^{\Delta}}, \nonumber \\
    \infnorm{\boldsymbol{h}_{k}^{(2), \Delta}} &\leq \alpha_{\delta}^{(2)} \infnorm{\boldsymbol{h}_{k-1}^{(2), \Delta}} + \beta_{\delta u}^{(2)} \infnorm{\boldsymbol{h}_{k}^{(1), \Delta}} \nonumber \\
    &\leq \alpha_{\delta}^{(2)} \infnorm{\boldsymbol{h}_{k-1}^{(2), \Delta}} + \alpha_{\delta}^{(1)} \beta_{\delta u}^{(2)} \infnorm{\boldsymbol{h}_{k-1}^{(1), \Delta}} + \nonumber \\
    &\quad + \beta_{\delta u}^{(2)} \beta_{\delta u}^{(1)} \infnorm{\boldsymbol{u}_{k-1}^{\Delta}}, \nonumber \\
    &\quad \vdots \nonumber \\
    \infnorm{\boldsymbol{h}_{k}^{(l), \Delta}} &\leq \alpha_{\delta}^{(l)} \infnorm{\boldsymbol{h}_{k-1}^{(l), \Delta}} + \sum_{\ell = 1}^{l-1} \alpha_{\delta}^{(\ell)} \! \prod_{\ell' = \ell + 1}^{l} \! \beta_{\delta u}^{(\ell')} \infnorm{\boldsymbol{h}_{k-1}^{(\ell), \Delta}} + \nonumber \\
    &\quad + \prod_{\ell = 1}^{l} \beta_{\delta u}^{(\ell)} \infnorm{\boldsymbol{u}_{k-1}^{\Delta}}.
\end{talign}
Now, define $\boldsymbol{\eta}_k^{\Delta} \defeq{} \left[\infnorm{\boldsymbol{h}_{k}^{(1), \Delta}}, \ldots, \infnorm{\boldsymbol{h}_{k}^{(L), \Delta}}\right]^\top \in \realset_{\geq 0}^{L}$ as the vector stacking all the norms of the $L$ hidden state differences at the time $k \in \naturalset$.
We can can combine all the inequalities in~\eqref{eq:incremental_hidden_state_upper_bound_all} to get:
\begin{talign}
    \label{eq:incremental_hidden_state_upper_bound_all_stacked}
    &\boldsymbol{\eta}_k^{\Delta} \leq A_{\delta} \boldsymbol{\eta}_{k-1}^{\Delta} + B_{\delta u} \infnorm{\boldsymbol{u}_{k-1}^{\Delta}},\\
    &A_\delta = \begin{bmatrix}
        \alpha_{\delta}^{(1)} & 0 & \cdots & 0 \\
        \alpha_{\delta}^{(1)} \beta_{\delta u}^{(2)} & \alpha_{\delta}^{(2)} & \cdots & 0 \\
        \vdots & \vdots & \ddots & \vdots \\
        \alpha_{\delta}^{(1)} \! \prod_{\ell = 2}^{L} \beta_{\delta u}^{(\ell)} & \alpha_{\delta}^{(2)} \! \prod_{\ell = 3}^{L} \beta_{\delta u}^{(\ell)} & \cdots & \alpha_{\delta}^{(L)}
    \end{bmatrix} \! \in \! \realset_{\geq 0}^{L \times L} \nonumber, \\
    &B_{\delta u} = \begin{bmatrix}
        \beta_{\delta u}^{(1)} & \beta_{\delta u}^{(2)} \beta_{\delta u}^{(1)} & \cdots & \prod_{\ell = 1}^{L} \beta_{\delta u}^{(\ell)}
    \end{bmatrix}^\top \! \in \! \realset_{\geq 0}^{L}. \nonumber
\end{talign}
Next, we propagate~\eqref{eq:incremental_hidden_state_upper_bound_all_stacked} backwards in a fashion similar to~\eqref{eq:hidden_state_upper_bound_propagated}~and~\eqref{eq:incrementa_hidden_state_upper_bound_propagated} to get (the $\leq$ operator is applied component-wise):
\begin{tequation}
    \label{eq:incremental_hidden_state_upper_bound_all_stacked_propagated}
    \boldsymbol{\eta}_k^{\Delta} \leq A_{\delta}^k \boldsymbol{\eta}_0^{\Delta} + \max_{0 \leq z < k} \infnorm{\boldsymbol{u}_{z}^{\Delta}} \sum_{z = 0}^{k - 1} A_{\delta}^{k - 1 - z} B_{\delta u}.
\end{tequation}
Now, notice that if the condition for layer-wise \dISS{} in~\eqref{eq:MGU_layer_dISS_condition} holds for every layer $l \in \layers$, we have $\alpha_{\delta}^{(l)} \in (0, 1)$, $\forall l \in \layers$.
Then, $A_\delta$ in~\eqref{eq:incremental_hidden_state_upper_bound_all_stacked}, being a lower triangular matrix, has all eigenvalues within the unit circle, making it Schur stable.
Moving on, consider the right-most summation in~\eqref{eq:incremental_hidden_state_upper_bound_all_stacked_propagated}.
Due to the non-negativeness of $A_{\delta}$ and $B_{\delta u}$, we have that (component-wise) $\sum_{z = 0}^{k - 1} A_{\delta}^{k - 1 - z} B_{\delta u} \leq \sum_{z = 0}^{\infty} A_{\delta}^{z} B_{\delta u}$.
However, $\sum_{z = 0}^{\infty} A_{\delta}^{z}$ is the Neumann series of a Schur stable matrix~\cite[Chapter 7]{meyer2023matrix}, which admits the closed-form solution $\sum_{z = 0}^{\infty} A_{\delta}^{z} = \left(I_L - A_{\delta}\right)^{-1}$, $I_L$ being the $L \times L$ identity matrix.
Then, going back to~\eqref{eq:incremental_hidden_state_upper_bound_all_stacked_propagated}, leveraging this finding and applying the norm to both sides of the inequality\footnote{Note that $\infnorm{\boldsymbol{h}_k^{\Delta}} = \infnorm{\left[\boldsymbol{h}_k^{(1), \Delta^\top}, \ldots, \boldsymbol{h}_k^{(L), \Delta^\top}\right]^\top} = \infnorm{\boldsymbol{\eta}_k^{\Delta}}$.}:%\footnote{Note that $\infnorm{\boldsymbol{h}_k^{\Delta}} = \infnorm{\left[\boldsymbol{h}_k^{(1), \Delta^\top}, \ldots, \boldsymbol{h}_k^{(L), \Delta^\top}\right]^\top} = \infnorm{\left[\infnorm{\boldsymbol{h}_{k}^{(1), \Delta}}, \ldots, \infnorm{\boldsymbol{h}_{k}^{(L), \Delta}}\right]^\top} = \infnorm{\boldsymbol{\eta}_k^{\Delta}}$.}:
\begin{talign}
    \label{eq:incremental_hidden_state_upper_bound_all_stacked_norm}
    \infnorm{\boldsymbol{h}_k^{\Delta}} &\leq \infnorm{A_{\delta}^k \boldsymbol{\eta}_0^{\Delta} + \max_{0 \leq z < k} \infnorm{\boldsymbol{u}_{z}^{\Delta}} \left(I_L - A_{\delta}\right)^{-1} B_{\delta u}} \nonumber \\
    &\leq \infnorm{A_{\delta}^k} \infnorm{\boldsymbol{h}_0^{\Delta}} + \infnorm{\left(I_L - A_{\delta}\right)^{-1} B_{\delta u}} \max_{0 \leq z < k} \infnorm{\boldsymbol{u}_{z}^{\Delta}}
\end{talign}
Concerning $\infnorm{A_{\delta}^k}$, making no assumption on the diagonalizability of $A_{\delta}$, we can leverage the Jordan form of $A_{\delta}$ to derive an upper bound for the aforementioned norm that is based on its spectral radius (see~\cite[Chapter 7]{meyer2023matrix})\footnote{If $A_{\delta}$ were to be diagonalizable, the inequality in~\eqref{eq:A_delta_k_norm_bound_1} would simplify to $\infnorm{A_{\delta}^k} \leq \omega \left(\max_{l\in\mathcal{L}}\alpha_{\delta}^{\left(l\right)}\right)^k$ for some $\omega \geq 1$.}:
\begin{tequation}
    \label{eq:A_delta_k_norm_bound_1}
    \infnorm{A_{\delta}^k} \leq \omega \sum_{z = 0}^{m - 1} \begin{pmatrix}k\\
z
\end{pmatrix}\left(\max_{l\in\mathcal{L}}\alpha_{\delta}^{\left(l\right)}\right)^{k-z}
\end{tequation}
for some $\omega \geq 1$, where $m \in \posintegerset$, $m \leq L$, is the size of the largest Jordan block of $A_{\delta}$, and, convention-wise, the binomial coefficient is such that $\begin{pmatrix}k\\
z
\end{pmatrix} = 0$ for any $z > k$.
Now, note that, for $k \to \infty$, the right-hand side of~\eqref{eq:A_delta_k_norm_bound_1} tends to $\omega \sum_{z = 0}^{m - 1} \frac{k^z}{z!} \left(\max_{l\in\mathcal{L}}\alpha_{\delta}^{\left(l\right)}\right)^{k - z}$.
Due to the Schur stability of $A_{\delta}$, it immediately follows that $\omega \sum_{z = 0}^{m - 1} \frac{k^z}{z!} \left(\max_{l\in\mathcal{L}}\alpha_{\delta}^{\left(l\right)}\right)^{k - z} \to 0$ for $k \to \infty$.
Thus, going back to~\eqref{eq:incremental_hidden_state_upper_bound_all_stacked_norm} and in the light of the previous considerations, we can conclude that
\begin{tequation}
    \label{eq:KL_dISS_network-wise}
    \psi_{\delta} \! \left(\!\infnorm{\boldsymbol{h}_{0}^{\mathrm{a}} \!\! - \! \boldsymbol{h}_{0}^{\mathrm{b}}} \! , k\!\right) \! \defeq{} \! \omega \! \sum_{z = 0}^{m - 1} \! \begin{pmatrix}k\\
z
\end{pmatrix} \! \left(\max_{l\in\mathcal{L}}\alpha_{\delta}^{\left(l\right)}\right)^{k-z} \!\infnorm{\boldsymbol{h}_{0}^{\mathrm{a}} \!\! - \! \boldsymbol{h}_{0}^{\mathrm{b}}} \!,
\end{tequation}
is the $\mathcal{KL}$ function for the \dISS{} condition in~\eqref{eq:diss}, while
\begin{tequation}
    \label{eq:Kinfty_dISS_network-wise}
    \gamma_{\delta u} \! \left(\!\max_{0 \leq z <k } \! \infnorm{\boldsymbol{u}_{z}^{\mathrm{a}} \!\! - \! \boldsymbol{u}_{z}^{\mathrm{b}} }\!\right) \! \defeq \! \infnorm{\left(I_L \! - \!A_{\delta}\right)^{-1} \! B_{\delta u}} \! \max_{0\leq z<k} \! \infnorm{\boldsymbol{u}_{z}^{\mathrm{a}} \! - \!\boldsymbol{u}_{z}^{\mathrm{b}}}
\end{tequation}
is the corresponding $\mathcal{K}_{\infty}$ function.
We can conclude that, if the condition for layer-wise \dISS{} in~\eqref{eq:MGU_layer_dISS_condition} holds for every layer $l \in \layers$, the resulting MGU network is \dISS{} in $\mathcal{H}_{\mathrm{inv}}$ with respect to $\mathcal{U}$ due to the existence of the functions~\eqref{eq:KL_dISS_network-wise}~and~\eqref{eq:Kinfty_dISS_network-wise}.
\hfill$\blacksquare$

%%%%%%%%%%%%%%%%%%%%%%%%%%%%%%%%%%%%%%%%%%%%%%%%%%%%%%%%%%%%%%%%%%%%%%%%%%%%%%%%%%%%%%%%%%%%%%%%%%%%%%%%%%

\paragraph*{Proof of Proposition~\ref{pr:MGU_network_dISS_imply_ISS}}
Theorem~\ref{th:MGU_network_dISS} requires every layer $l \in \layers$ of the MGU network to satisfy the \dISS{} condition in~\eqref{eq:MGU_layer_dISS_condition}.
Proposition~\ref{pr:MGU_layer_dISS_imply_ISS} establishes that, for any layer $l$, satisfaction of~\eqref{eq:MGU_layer_dISS_condition} implies satisfaction of the ISS condition in~\eqref{eq:MGU_layer_ISS_condition}.
Therefore, all layers satisfy the conditions of Theorem~\ref{th:MGU_network_ISS}, guaranteeing the MGU network to be ISS.
\hfill$\blacksquare$

%%%%%%%%%%%%%%%%%%%%%%%%%%%%%%%%%%%%%%%%%%%%%%%%%%%%%%%%%%%%%%%%%%%%%%%%%%%%%%%%%%%%%%%%%%%%%%%%%%%%%%%%%%

\paragraph*{Proof of Proposition~\ref{pr:convex_projection_implies_iss}}
Consider all $l \in \layers$ and $\feasibleset$ in~\eqref{eq:iss_feasible_set}.
For any $\parameters \in \feasibleset$, we have $\infnorm{R_{\tilde{h}}^{(l)}} \leq 1 - \projectionmargin$.
Since $\bar{\sigma}_{f}^{(l)} \in (0,1)$ due to the range of the sigmoid function and $\projectionmargin \in (0, 1)$, we have:
\begin{equation*}
    \bar{\sigma}_{f}^{(l)} \infnorm{R_{\tilde{h}}^{(l)}} < \infnorm{R_{\tilde{h}}^{(l)}} \leq 1 - \projectionmargin < 1.
\end{equation*}
Consequently, any $\parameters \in \feasibleset$ satisfies the layer-wise ISS condition in~\eqref{eq:MGU_layer_ISS_condition} for every $l \in \layers$. 
Thus, due to Theorem~\ref{th:MGU_network_ISS}, the resulting MGU network is ISS.
\hfill$\blacksquare$

\bibliographystyle{plain}
\bibliography{bibliography/bibliography}
\end{document}